%
\font\sc=cmcsc10
\font\tencmmib=cmmib10
\font\eightcmmib=cmmib10
\font\tenmsym=msbm10	
\font\eightmsym=msbm8	
\font\teneufm=eufm10
\font\tencmsy=cmsy10
\font\eightcmsy=cmsy8
\textfont 7=\tencmsy \scriptfont 7=\eightcmsy 
\def\scr {\fam7 } 
\textfont 9=\tenmsym \scriptfont 9=\eightmsym 
\def\bb {\fam9 } 
\textfont 8=\teneufm \scriptfont 8=\teneufm 
\def\frak {\fam8 } 
\textfont 6=\tencmmib \scriptfont 7=\eightcmmib 
\def\bi {\fam6 } 
\magnification=1200
\font\bgf=cmbx10 scaled\magstep2

\font\tt=cmtt10
\def\bib#1{\noindent\hbox to\parindent{[#1]\hfil}\hang}

\centerline{\bgf Construction of Miniversal Deformations of Lie Algebras}
\bigskip

\centerline{\sc Alice FIALOWSKI \qquad and \qquad Dmitry FUCHS} \medskip

\noindent
\hbox{\hskip50pt
\vtop{
   \hbox to 2in {\hfil Department of Applied Analysis\hfil}\par
   \hbox to 2in {\hfil E\"otv\"os Lor\'and University\hfil} \par
   \hbox to 2in {\hfil M\'uzeum krt.~6-8\hfil}\par
   \hbox to 2in {\hfil H-1088 Budapest\hfil} \par
   \hbox to 2in {\hfil Hungary\hfil}}
\quad
\vtop{
   \hbox to 2in{\hfil Department of Mathematics\hfil} \par
   \hbox to 2in {\hfil University of California\hfil} \par
   \hbox to 2in {\hfil Davis CA 95616\hfil} \par
   \hbox to 2in {\hfil USA\hfil}}
}

\vskip 25pt

\centerline{\bf Introduction}

\medskip

In this paper we consider deformations of finite or infinite
dimensional Lie algebras over a field of characteristic
$0$. By ``deformations of a Lie algebra'' we mean the (affine
algebraic) manifold of all Lie brackets.  Consider the quotient of
this variety by the action of the group GL.  It is well-known (see
[Hart]) that in the category of algebraic varieties the quotient by a
group action does not always exist.  Specifically, there is in general
no universal deformation of a Lie algebra $L$ with a commutative
algebra base $A$ with the property that for any other deformation of
$L$ with base $B$ there exists a unique homomorphism $f\colon B\to A$
that induces an equivalent deformation.  If such a homomorphism
exists (but not unique), we call the deformation of $L$ with base $A$
{\it versal}.

\smallskip

Classical deformation theory of associative and Lie algebras began
with the works of Gerstenhaber [G] and Nijenhuis-Richardson [NR] in
the 1960's. They studied one-parameter deformations and established
the connection between Lie algebra cohomology and infinitesimal
deformations.  They did not study the versal property of deformations.

\smallskip

A more general deformation theory for Lie algebras follows from
Schlessinger's work [Sch].  If we consider deformations with base
spec$\,A$, where $A$ is a local algebra, this set-up is adequate to
study the problem of ``universality'' among formal deformations.  This
was worked out for Lie algebras in [Fi1], [Fi3]; it turns out that in
this case under some minor restrictions there exists a so-called {\it
miniversal} element.  The problem is to construct this element.

\smallskip

There is confusion in the literature when one tries to describe {\it
all} the nonequivalent deformations of a given Lie algebra.  There
were several attempts to work out an appropriate theory for solving this
basic problem in deformation theory, but none of them were completely adequate.

\smallskip

The construction below is parallel to the general constructions in
deformation theory, as in [P], [I], [La], [GoM], [K]. The general
theory, which can provide a construction of a local miniversal
deformation, is outlined in [Fi1]. The procedure however needs a
proper theory of Massey operations in the cohomology, and an algorithm
for computing all the possible ways for a given infinitesimal
deformation to extend to a formal deformation.  The proper theory of
Massey operations is developed in [FuL].  Our understanding of the
construction arose from the study of the infinite dimensional Lie
algebra $L_1$ of polynomial vector fields in $\bb C$ with trivial
1-jet at $0$, in which case we completely described a miniversal
deformation. In [FiFu] we proved that the base of the miniversal
deformation of this Lie algebra is the union of three algebraic
curves: two smooth curves and another curve with a cusp at $0$, with
the tangent lines to all three curves coinciding at $0$. \medskip

The structure of the paper is as follows: In Section 1 we give the
necessary definitions and some facts on infinitesimal deformations. In
Section 2 we recall Harrison cohomology, in Section 3 discuss
obstruction theory. Section 4 gives the theoretical construction of a
miniversal deformation, and some preliminary computations. Section 5
recalls the proper Massey product definition and describes its
properties (see [FuL]).  In Section 6 we calculate
obstructions. Section 7 provides a scheme for computing the base of a
miniversal deformation of a Lie algebra convenient for practical
use. In Section 8 we apply the construction to the Lie algebra $L_1$. \medskip

\noindent{\sc Acknowledgements}

We thank the Erwin Schr\"odinger Institute for Mathematical Physics,
Vienna, where the work began, the first author thanks the
Max-Planck-Institute f\"ur Mathematik, Bonn, for hospitality and partial
support. Both authors thank Michael Penkava for reading the manuscript,
correcting a number of misprints and making valuable suggestions.
\bigbreak

\centerline{\bf 1. Lie algebra deformations} \medskip

{\bf 1.1.} Let $L$ be a Lie algebra over a characteristic 0 field $\bb K$, 
and let $A$ be a commutative algebra with identity over $\bb K$ with a fixed 
augmentation $\varepsilon\colon A\to{\bb K},\ \varepsilon(1)=1$; we set ${\rm 
Ker}\, \varepsilon=\frak m$. To avoid transfinite induction, we will
 assume that $\dim({\frak m}^k/{\frak m}^{k+1})<\infty$ for all $k$.
\smallskip 

{\sc Definition} 1.1. A {\it deformation} $\lambda$ of $L$ with base 
$(A,{\frak m})$, or simply with base $A$, is a Lie $A$-algebra 
structure  on the tensor product
$A\otimes_{\bb K}L$ with the bracket $[\, ,\, ]_\lambda$, such that 
$$\varepsilon\otimes{\rm id}\colon A\otimes L\to {\bb K}\otimes L=L$$ 
is a Lie algebra homomorphism. (We usually abbreviate $\otimes_{\bb 
K}$ to $\otimes$.) See [Fi1], [Fi3].
\smallskip 

{\sc Example} 1.2. If $A={\bb K}[t]$, then a deformation of $L$ with 
 base $A$ is the same as an algebraic 1-parameter deformation of 
$L$. More generally, if $A$ is the algebra of regular functions on an 
affine algebraic manifold $X$, then a deformation of $L$ with base 
$A$ is the same as an algebraic deformation of $L$ with base 
$X$.\smallskip 

Two deformations of a Lie algebra $L$ with the same base $A$ are 
called {\it equivalent} (or isomorphic) if there exists a Lie algebra 
isomorphism between the two copies of $A\otimes L$ with the two Lie 
algebra structures, compatible with $\varepsilon\otimes{\rm id}$. A 
deformation with base $A$ is called {\it local} if the algebra $A$ 
is local, and it is called {\it infinitesimal} if, in addition to 
this, ${\frak m}^2=0$.\smallskip 

{\sc Definition} 1.3. Let $A$ be a complete local algebra (completeness
means that $A=\overleftarrow{\displaystyle{\lim_{n\to\infty}}}(A/{\frak m}
^n)$, where $\frak m$ is the maximal ideal in $A$). A {\it formal 
deformation} of $L$ with  
base $A$ is a Lie $A$-algebra structure on the completed tensor product $A\, 
\hat\otimes\,  L=\overleftarrow{\displaystyle{\lim_{n\to\infty}}}\, ((A/{\frak 
m}^n)\otimes L)$ such that $$\varepsilon\, \hat\otimes\, {\rm id}\colon A\, 
\hat\otimes\,  L\to {\bb K}\otimes L=L$$is a Lie algebra homomorphism
(see [Fi3]) . 

The above notion of equivalence is extended to formal deformations in an 
obvious way.\smallskip 

{\sc Example} 1.4. If $A={\bb K}[[t]]$ then a formal deformation of 
$L$ with base $A$ is the same as a formal 1-parameter deformation 
of $L$. See [G], [NR].
\bigskip 

Let $A'$ be another commutative algebra with identity over $\bb K$ with a fixed 
augmentation $\varepsilon'\colon A'\to\bb K$, and let $\varphi\colon A\to A'$ 
be an algebra homomorphism with $\varphi(1)=1$ and 
$\varepsilon'\circ\varphi=\varepsilon$.\smallskip 

{\sc Definition} 1.5. If a deformation $\lambda$ of $L$ with base 
$(A,{\frak m})$ is given, then the {\it push-out} 
$\varphi_\ast\lambda$ is the deformation of $L$ with base 
$(A',{\frak m}'={\rm Ker}\, \varepsilon')$, which is the Lie algebra 
structure $$[a'_1\otimes_A(a_1\otimes l_1),a'_2\otimes_A(a_2\otimes 
l_2)]'=a'_1a'_2\otimes_A[a_1\otimes l_1,a_2\otimes 
l_2],$$$$a'_1,a'_2\in A',\ a_1,a_2\in A,\ l_1,l_2\in L,$$on $A'\otimes 
L=(A'\otimes_AA)\otimes L=A'\otimes_A(A\otimes L)$. Here $A'$ is 
regarded as an $A$-module with the structure $a'a=a'\varphi(a)$, and 
the operation $[\, ,\, ]$ in the right hand side of the formula refers 
to the Lie algebra structure $\lambda$ on $A\otimes L$. 

The push-out of a formal deformation is defined in a similar way.
\medskip

{\bf 1.2.} For completeness' sake, we recall the definition of Lie algebra 
cohomology (see [Fu]). We need only the case of cohomology with 
coefficients in the adjoint representation, and therefore we restrict 
our definition to this case. 

Let$$C^q(L;L)={\rm Hom}\, (\Lambda^qL,L)$$be the space of all skew-symmetric 
$q$-linear forms on a Lie algebra $L$ with values in $L$. Define the
differential$$\delta\colon C^q(L;L)\to C^{q+1}(L;L)$$by the formula
$$\eqalign{(\delta\gamma)(l_1,\dots ,l_{q+1})&=\sum_{1\le s<t\le q+1}
(-1)^{s+t-1}\gamma([l_s,l_t],l_1,\dots\hat l_s\dots\hat l_t\dots ,l_{q+1})\cr
&+\sum_{1\le u\le q+1}[l_u,\gamma(l_1,\dots\hat l_u\dots,l_{q+1})],\cr}$$
where $\gamma\in C^q(L;L),\ l_1,\dots ,l_{q+1}\in L$. It may be checked that
$\delta^2=0$, and the cohomology of the complex $\{C^q(L;L),\delta\}$ 
is denoted by $H^q(L;L)$. 

For $\alpha\in C^p(L;L),\, \beta\in C^q(L;L),$ the Lie product $[\alpha,\beta]$
is defined be the formula$$\eqalign{[\alpha,\beta]&(l_1,\dots ,l_{p+q-1})\cr
&=\sum_{1\le j_1<\dots<j_q\le p+q-1}(-1)^{\Sigma_s(j_s-s)}\alpha(\beta(l_{j_1},
\dots,l_{j_q}),l_1,\dots\hat l_{j_1}\dots\hat l_{j_q}\dots,l_{p+q-1})-\cr}$$$$
(-1)^{(p-1)(q-1)}\sum_{1\le k_1<\dots<k_p\le p+q-1}(-1)^{\Sigma_t(k_t-t)}\beta
(\alpha(l_{k_1},\dots,l_{k_p}),l_1,\dots\hat l_{k_1}\dots\hat l_{k_p}\dots,
l_{p+q-1}).$$ If one sets ${\scr C}^q=C^{q+1}(L;L),\, {\scr H}^q=H^{q+1}(L;L)$, 
then this bracket operation (with the differential $\delta$) makes ${\scr C}=
\bigoplus{\scr C}^q$ a differential graded Lie algebra (DGLA), and makes 
${\scr H}=\bigoplus{\scr H}^q$ a graded Lie algebra. \medskip

{\bf 1.3.} Here is the fundamental example of an infinitesimal 
deformation of a Lie algebra. Consider a Lie algebra $L$ which 
satisfies the condition$$\dim H^2(L;L)<\infty .$$ This is true, for 
example, if $\dim L<\infty$. 

(There are some ways to weaken if not to completely avoid this condition. 
For example, if the Lie algebra $L$ is $\bb Z$-graded, $L=\bigoplus_{q\in\bb 
Z}L_{(q)},\ [L_{(p)},L_{(q)}]\subset L_{(p+q)},$ then $H^2(L;L)$ also becomes 
graded, $H^2(L;L)=\bigoplus_{q\in\bb Z}H^2_{(q)}(L;L)$, and the construction will 
be valid in a slightly modified form, if one supposes that $\dim 
H^2_{(q)}(L;L)<\infty$ for all $q$. See the details in 7.4 below.)

Consider the algebra$$A={\bb K}\oplus H^2(L;L)'$$ with the second summand 
being an ideal with zero multiplication ($'$ means the dual). Fix some 
homomorphism$$\mu\colon H^2(L;L)\to C^2(L;L)={\rm Hom}\, (\Lambda^2L,L)$$ which 
takes a cohomology class into a cocycle representing this class. Define a 
Lie algebra structure on $$A\otimes L=({\bb K}\otimes L)\oplus(H^2(L;L)'\otimes 
L)=L\oplus{\rm Hom}\, (H^2(L;L),L)$$ by the 
formula$$[(l_1,\varphi_1),(l_2,\varphi_2)]=([l_1,l_2],\psi),$$ where 
$$\psi(\alpha)=\mu(\alpha)(l_1,l_2)+[\varphi_1(\alpha),l_2]+[l_1,\varphi_2 
(\alpha)],$$$$l_1,l_2\in L,\ \varphi_1,\varphi_2\in{\rm Hom}\, 
(H^2(L;L),L)),\, \alpha\in H^2(L;L).$$ (The Jacobi identity for this operation 
is implied by $\delta\mu(\alpha)=0$.) This determines a deformation of $L$ 
with base $A$ which is clearly infinitesimal.\smallskip 

{\sc Proposition} 1.6. {\it Up to an isomorphism, this deformation 
does not depend on the choice of $\mu$.}\smallskip 

{\sc Proof}. Let $$\mu'\colon H^2(L;L)\to C^2(L;L)$$ be another choice
for $\mu$. Then there exists a homomorphism$$\gamma\colon H^2(L;L)\to
C^1(L;L)={\rm Hom}\, (L,L)$$such that
$\mu'(\alpha)=\mu(\alpha)+\delta\gamma(\alpha)$ for all $\alpha\in
H^2(L;L)$. Define a linear automorphism $\rho$ of the space $A\otimes
L=L\oplus{\rm Hom}\, (H^2(L;L),L)$ by the
formula$$\rho(l,\varphi)=(l,\psi),\
\psi(\alpha)=\varphi(\alpha)+\gamma(\alpha)(l),$$$$l\in L,\ \varphi\in{\rm
Hom}\, (H^2(L;L),L),\ \alpha\in H^2(L;L).$$ The map $\rho$ is clearly an
automorphism. The inverse of $\rho$ is given by replacing $\gamma$ with
$-\gamma$ in the formula. To prove that $\rho$ is an isomorphism
between the two Lie algebra structures, one needs to check that for any
$l_1,l_2\in L,\, \varphi_1,\varphi_2\in{\rm Hom}\, (H^2(L;L),L),\,
\alpha\in H^2(L;L)$ one has
$$\eqalign{\mu(\alpha)(l_1,l_2)&+[\varphi_1(\alpha),l_2]+[l_1,
\varphi_2(\alpha)]+\gamma(\alpha)([l_1,l_1])\cr
&=\mu'(\alpha)(l_1,l_2)+[\varphi_1(\alpha)+\gamma(\alpha)(l_1),l_2]+[l_1,
\varphi_2(\alpha)+\gamma(\alpha)(l_2)];\cr}$$but this follows directly
from the equality
$\mu'(\alpha)=\mu(\alpha)+\delta\gamma(\alpha)$.\smallskip

We will denote the infinitesimal deformation of $L$ constructed above by
$\eta_L$.  \medskip

{\bf 1.4.} The main property of $\eta_L$ is its (co-)universality in the 
class of infinitesimal deformations.\smallskip

Let $\lambda$ be an infinitesimal deformation of the Lie algebra $L$
with the finite dimensional base $A$. Take $\xi\in{\frak m}'$, or,
equivalently, $\xi\in A'$ and $\xi(1)=0$. For $\l_1,l_2\in L$
set 
$$\alpha_{\lambda,\xi}(l_1,l_2)=(\xi\otimes{\rm id})[1\otimes
l_1,1\otimes l_2]_\lambda\in{\bb K}\otimes L=L.$$

{\sc Lemma 1.7}. {\it The cochain $\alpha_{\lambda,\xi} \in C^2(L;L)$ is
a cocycle}.  \medskip

{\sc Proof.} Let $l_1, l_2, l_3 \in L$. Since $
[1\otimes l_1,1\otimes l_2]_\lambda - 1\otimes [l_1,l_2] \in {\frak
m}\otimes L$,
we have
$$
[1\otimes l_1,1\otimes l_2]_\lambda = 1\otimes [l_1,l_2] + \sum_i m_i
\otimes k_i,
$$
where $m_i \in {\frak m}, k_i \in L$. Hence
$$
(\xi \otimes {\rm id}) [1\otimes l_1, 1\otimes l_2]_\lambda = (\xi \otimes {\rm id})
[1\otimes [l_1,l_2], 1\otimes l_3]_\lambda + 
(\xi \otimes {\rm id})\sum_i m_i[1\otimes k_i, 1 \otimes l_3].
$$
The first summand here is $\alpha_{\lambda,\xi}([l_1,l_2],l_3)$.  For
the second summand
$$
m_i[1\otimes k_i,1\otimes l_3] = m_i(1\otimes [k_i,l_3] + h),
$$
where $h\in {\frak m}\otimes L$.  Since ${\frak m}^2=0$ we have $m_ih=0$.
Hence $m_i[1\otimes k_i, 1\otimes l_3] = m_i\otimes [k_i,l_3]$,
and
$$\eqalign{
(\xi \otimes {\rm id})\sum_i m_i[1\otimes k_i, 1 \otimes l_3]
&= \sum_i(\xi \otimes {\rm id})(m_i\otimes[k_i,l_3]) 
	=\sum_i\xi(m_i)[k_i,l_3] \cr
&= \sum_i [\xi(m_i)k_i,l_3] 
	=[(\xi\otimes {\rm id})\Bigl(\sum_i m_i\otimes k_i\Bigr), l_3] \cr
&= [(\xi \otimes {\rm id})([1\otimes l_1, 1\otimes l_2]_\lambda -
	1\otimes[l_1, l_2]), l_3] \cr
&= [(\xi \otimes {\rm id})[1\otimes l_1, 1\otimes l_2], l_3] 
= [\alpha_{\lambda,\xi}(l_1, l_2), l_3].}
$$
In the last step above we used that $\xi(1) = 0$.  Thus
$$
(\xi\otimes {\rm id})[[1\otimes l_1, 1\otimes l_2]_\lambda, 1\otimes
l_3]_\lambda
 = \alpha_{\lambda,\xi}([l_1,l_2],l_3) + [\alpha_{\lambda,\xi}(l_1,
l_2), l_3],
$$
and the Jacobi identity for $[\ , \ ]_\lambda$ shows that
$\alpha_{\lambda,\xi}$ is a cocycle.  \medskip

{\sc Proposition 1.8}. {\it For any infinitesimal deformation $\lambda$
of the Lie algebra $L$ with a finite-dimensional base $A$
there exists a unique homomorphism $\varphi\colon{\bb K}\oplus
H^2(L;L)'\to A$ such that $\lambda$ is equivalent to the push-out
$\varphi_\ast\eta_L$.} 
\smallskip

{\sc Proof.} For $\xi\in{\frak m}'$ let $a_{\lambda,\xi}\in H^2(L;L)$ be
the cohomology class of the cocycle $\alpha_{\lambda,\xi}$. The
correspondences $\xi\mapsto\alpha_{\lambda,\xi}, \xi\mapsto
a_{\lambda,\xi}$ define homomorphisms $$\alpha_\lambda\colon{\frak
m}'\to C^2(L;L),\qquad \delta\circ\alpha_\lambda=0,$$
$$a_\lambda\colon{\frak m}'\to H^2(L;L).$$We claim that

(i) the deformation $\lambda$ is fully determined by $\alpha_\lambda$;

(ii) the deformations $\lambda,\lambda'$ are equivalent if and only if 
$a_\lambda=a_{\lambda'}$;

(iii) if $\varphi={\rm id}\oplus a'_\lambda\colon{\bb K}\oplus 
H^2(L;L)'\to{\bb K}\oplus{\frak m}=L$, then $\varphi_\ast\eta_L$ is
equivalent to $\lambda$.

\noindent Since (ii) and (iii) obviously imply Proposition, it remains to 
prove (i)--(iii). The statement (i) is obvious. To prove (ii) notice that 
an $A$-automorphism $\rho\colon A\otimes L\to A\otimes L$, that is
$$L\oplus({\frak m}\otimes L)\to L\oplus({\frak m}\otimes L),$$whose
$L\to L$ part is the identity (this is the condition of
compatibility with  $\varepsilon\otimes{\rm id}$), is fully
determined by its $$L\to{\frak  m}\otimes L$$part, which we denote
by $b_\rho$, and the latter may  be chosen arbitrarily. This is an
element of$${\rm Hom}\, (L,{\frak  m}\otimes L)={\frak m}\otimes{\rm
Hom}\, (L,L)={\frak m}\otimes C^1(L;L)= {\rm Hom}\, ({\frak
m}',C^1(L;L)).$$It is easy to check that $\rho$ establishes an
isomorphism between the Lie algebra structures $\lambda$ and 
$\lambda'$ if and only if
$$
\alpha_{\lambda'}-\alpha_\lambda=\delta\circ b_\rho,
$$
which proves (ii). Finally, it follows from the definitions that
$$\alpha_{\varphi_\ast\eta_L}=\mu\circ a_\lambda,$$which implies that
$a_{\varphi_\ast\eta_L}= a_\lambda$, and hence $\varphi_\ast\eta_L$
and $\lambda$ are isomorphic as was stated in (iii).\smallskip

{\sc Remark} 1.9. Technically, the mapping $a_\lambda\colon{\frak m}'\to 
H^2(L;L)$ constructed in the proof will be more important for us than 
the map $\varphi={\rm id}\oplus a'_\lambda$.\smallskip

Let $A$ be a local algebra with ${\rm dim}\, (A/{\frak m}^2)<\infty$.
Obviously, $A/{\frak m}^2$ is local with the maximal ideal ${\frak
m}/{\frak m}^2$, and $({\frak m}/{\frak m}^2)^2=0$. Recall that the
dual space $({\frak m}/{\frak m}^2)'$ is called the {\it tangent
space} of $A$; we denote it by $TA$. \smallskip

{\sc Definition} 1.10. Let $\lambda$ be a deformation of $L$ with base 
$A$. Then the mapping$$a_{\pi_\ast\lambda}\colon TA=({\frak m}/{\frak m}^2)'
\to H^2(L;L),$$where $\pi$ is the projection $A\to A/{\frak m}^2$, is called
the {\it differential} of $\lambda$ and is denoted by $d\lambda$.

The differential of a formal deformation is defined in a similar way.
\smallskip

It is clear from the construction that equivalent deformations or formal 
deformations have equal differentials. \medskip

{\bf 1.5.} It is not possible to construct a local or formal deformation 
of a Lie algebra with a similar universality property in the class of local 
or formal deformations. But it becomes possible for an appropriate 
weakening of this property.\smallskip

{\sc Definition} 1.11. A formal deformation $\eta$ of a Lie algebra $L$ 
with base $B$ is called {\it miniversal} if

(i) for any formal deformation $\lambda$ of $L$ with any (local) base $A$
there exists a homomorphism $f\colon B\to A$ such that the deformation 
$\lambda$ is equivalent to $f_\ast\eta$;

(ii) in the notations of (i), if $A$ satisfies the condition ${\frak m}^2=0$, 
then $f$ is unique (see [Fi1]).

If $\eta$ satisfies only the condition (i), then it is called {\it versal}.
\smallskip

Our goal is to construct a miniversal formal deformation of a given Lie 
algebra.\bigbreak

\centerline{\bf 2. Harrison cohomology.} \medskip

{\bf 2.1.} We will need a special cohomology theory for commutative 
algebras introduced in 1961 by D.~K.~Harrison [Harr]. The following general 
definition is contained in the article [B].

Let $A$ be a commutative $\bb K$-algebra. Consider the standard Hochschild 
complex $\{C_q(A),\partial\}$ for $A$. Here $C_q(A)$ is the $A$-module 
$A^{q+1}=A\otimes\dots\otimes A$ ($q+1$ factors), $A$ operates on the last 
factor, and the differential $\partial\colon C_q(A)\to C_{q-1}(A)$ is defined
by the formula$$\eqalign{\partial[a_1,\dots,a_q]=a_1[a_2,\dots,a_q]&+
\sum_{i=1}^{q-1}(-1)^q[a_1,\dots,a_ia_{i+1},\dots,a_q]\cr 
&+(-1)^qa_q[a_1,\dots,a_{q-1}],\cr}$$where $b_0[b_1,\dots,b_n]$ means
$b_0\otimes b_1\otimes\dots\otimes b_n\in C_n(A)$. A permutation from $S(q)$
is called a $(p,q-p)$-shuffle if the inverse permutation $(j_1,\dots,j_q)$
satisfies the conditions $j_1<\dots<j_p,\, j_{p+1}<\dots<j_q$. Let ${\rm Sh}
\, (p,q-p)\subset S(q)$ be the set of all $(p,q-p)$-shuffles. For $a_1,\dots,
a_q\in A$ and $0<p<q$ set$$s_p(a_1,\dots,a_q)=\sum_{(i_1,\dots,i_q)\in
{\rm Sh}(p,q-p)}{\rm sgn}(i_1,\dots,i_q)[a_{i_1},\dots,a_{i_q}]\in C_q(A).$$
Let ${\rm Sh}_q(A)$ be the $A$-submodule of $C_q(A)$ generated by 
the chains $s_p(a_1, \dots,a_q)$ for all $a_1,\dots,a_q\in A,\, 
0<p<q$. It may be checked (see [B], Proposition 2.2) that 
$\partial({\rm Sh}_q(A))\subset{\rm Sh}_{q-1}(A),$ which yields a 
complex ${\rm Ch}\, (A)=\{{\rm Ch}_q(A)=C_q(A)/{\rm Sh}_q(A), 
\partial\}$. This is the {\it Harrison complex}.\smallskip 

{\sc Definition} 2.1. For an $A$-module $M$ we set$$\eqalign{H_q^{\rm Harr}
(A;M)&=H_q({\rm Ch}\, (A)\otimes M),\cr H^q_{\rm Harr}(A;M)&=H^q({\rm Hom}\, 
({\rm Ch}\, (A),M);\cr}$$these are {\it Harrison homology and cohomology of
$A$ with coefficients in $M$.} (For the relations between Harrison and 
Hochschild homology and cohomology see [B].)\smallskip

We will need the following standard fact, which follows directly from 
the definition.\smallskip 

{\sc Proposition} 2.2. {\it Let $A$ be a local commutative $\bb K$-
algebra with the maximal ideal $\frak m$, and let $M$ be an $A$-module 
with ${\frak m}M=0$. Then we have the canonical isomorphisms 
$$H_q^{\rm Harr}(A;M)\cong H_q^{\rm Harr}(A;{\bb K})\otimes M,\ H^q_{\rm 
Harr}(A;M)\cong H^q_{\rm Harr}(A;{\bb K})\otimes M.$$} 

{\bf 2.2.} We will need only 1- and 2-dimensional Harrison cohomology. Here
is their direct description (belonging to Harrison [Harr]). Let $A$ and $M$ be
as above. Consider the complex$$0\to {\rm Ch}^1\, {\buildrel 
d_1\over\longrightarrow}\, {\rm Ch}^2\, {\buildrel
d_2\over\longrightarrow}\, C^3,$$where$${\rm Ch}^1={\rm Hom}\, (A,M),\ {\rm
Ch}^2={\rm Hom}\, (S^2A,M),\ C^3={\rm Hom}\, (A\otimes A\otimes A,M),$$$$
d_1\psi(a,b)=a\psi(b)-\psi(ab)+b\psi(a),\ d_2\varphi(a,b,c)=
a\varphi(b,c)-\varphi(ab,c)+\varphi(a,bc)-c\varphi(a,b),$$$$m\in M,\ a,b,c
\in A,\ \psi\in{\rm Ch}^1,\ \varphi\in{\rm Ch}^2.$$

{\sc Proposition} 2.3. (i) {\it $H^1_{\rm Harr}(A;M)$ is the space of 
derivations $A\to M$.} (ii) {\it Elements of $H^2_{\rm Harr}(A;M)$  
correspond bijectively 
to isomorphism classes of extensions $0\to M\to B\to A\to 0$ of the 
algebra $A$ by means of $M$.}\smallskip

{\sc Proof}. Part (i) is obvious. To prove (ii), consider an extension
$0\to M\, {\buildrel i\over\longrightarrow}\, B\, {\buildrel p\over 
\longrightarrow}\, A\to 0$ and fix a section $q\colon A\to B$ of $p$. 
Then $b\mapsto(p(b),i^{-1}(b-q\circ p(b)))$ is an isomorphism $B\to 
A\oplus M$. Let $(a,m)_q\in B$ be the inverse image of $(a,m)\in A\oplus M$ 
with respect to this isomorphism. For $a_1,a_2\in A$ set $\varphi_q(a_1,a_2)
=i^{-1}((a_1,0)_q(a_2,0)_q-(a_1a_2,0)_q)\in M$. Then the multiplication in
$B$ is $(a_1,m_1)_q(a_2,m_2)_q=(a_1a_2,a_1m_2+a_2m_1+\varphi_q(a_1,a_2))_q$,
so it is determined by $\varphi_q$. Furthermore, the associativity of 
the algebra $B$ implies that $\varphi_q\in{\rm Ch}^2$ is a cocycle. 
For any other section $q'\colon A\to B$ one has $i^{-1}\circ(q'-q) \in 
{\rm Ch}^1$, and it is easy to check that 
$\varphi_{q'}=\varphi_q+d_1(i^{-1}\circ(q'-q))$. This implies 
(ii).\smallskip 

{\sc Corollary} 2.4. {\it If $A$ is a local algebra with the maximal ideal 
$\frak m$, then $H^1_{\rm Harr}(A;{\bb K})=({\frak m}/{\frak m}^2)'=TA$.}
\smallskip

{\sc Proof} Let $\varphi\colon A\to{\bb K}$ be a derivation. If $a\in{\frak
m}^2$, that is $a=a_1a_2,\ a_1,a_2\in{\frak m}$, then $\varphi(a)=
\varphi(a_1a_2)=a_1\varphi(a_2)+a_2\varphi(a_1)=0$ (since ${\frak m}{\bb K}
=0$). Furthermore, $\varphi(1)=\varphi(1\cdot1)=1\varphi(1)+1\varphi(1)=
2\varphi(1)$, hence $\varphi(1)=0$. On the other hand, any homomorphism 
$\varphi\colon A\to\bb K$ such that $\varphi({\frak m}^2)=0,\varphi(1)=0$, 
is a derivation. Thus the space of derivations $A\to\bb K$ is $({\frak m}/
{\frak m}^2)'$.\smallskip

{\sc Proposition} 2.5. {\it Let $0\to M\, {\buildrel
i\over\longrightarrow}\, B\, {\buildrel p\over\longrightarrow}\, A\to 0$
be an extension of an algebra $A$ .} (i) {\it If $A$ has an identity, then so
does $B$.} (ii) {\it If $A$ is local with the maximal ideal $\frak m$,
then $B$ is local with the maximal ideal $p^{-1}({\frak m})$.}\smallskip

{\sc Proof.} (i) We use the notations of the previous proof. Fix a section 
$q\colon A\to B$ of $p$. Then we get a cocycle $\varphi=\varphi_q\in{\rm
Ch}^2$. For any $a\in A$$$d_2\varphi(1,1,a)=\varphi(1,a)-\varphi(1,a)+
\varphi(1,a)-a\varphi(1,1)=0,$$which shows that $\varphi(1,a)=a\varphi(1,1)$.
Consider an arbitrary $\psi\in{\rm Ch}^1$ with $\psi(1)=\varphi(1,1)$. Let
$\varphi'=\varphi-d_1\psi$. Then for any $a\in A$
$$\eqalign{\varphi'(1,a)&=
\varphi(1,a)-d_1\psi(1,a)\cr &=\varphi(1,a)-\psi(a)+\psi(a)-a\psi(1)\cr
&=\varphi(1,a)-a\varphi(1,1)=0.\cr}$$
According to the previous proof, 
$\varphi'=\varphi_{q'}$ for some section $q'\colon A\to B$, and one has
$$(1,0)_{q'}(a,m)_{q'}=(a,m+\varphi_{q'}(1,a))_{q'}=(a,m)_{q'}.$$Hence,
$(1,0)_{q'}\in B$ is the unit element.

(ii) Let ${\frak n}\subset B$ be an ideal, and let ${\frak n}\not\subset
p^{-1}({\frak m})$. Then there is some $b\in\frak n$ such that $p(b)=1$. 
Choose a section $q\colon A\to B$ with $q(1)=b.$ Then $b=(1,0)_q$. For any
$(a,m)_q\in B$ one has$$(a,m)_q=(1,0)_q(a,m-\varphi_q(1,a))\in{\frak n},$$
and hence ${\frak n}=B$. \medskip

{\bf 2.3.} The relationship between the second Harrison cohomology of a
finite-dimensional local commutative algebra $A$ and extensions of $A$ 
may be also described in terms of one remarkable extension. This is the 
extension $$0\to H^2_{\rm Harr}(A;{\bb K})'\to C\to A\to 0,\eqno(1)$$ where 
the operation of $A$ on $H^2_{\rm Harr}(A;{\bb K})'$ is induced by the 
operation of $A$ on $\bb K$, and the cocycle 
$$
    f_A: S^2A\to H^2_{\rm Harr}(A;{\bb K})'
$$
is defined as the dual of a homomorphism $$\mu: H^2_{\rm Harr}
(A;{\bb K})\to {\rm Ch}^2(A;{\bb K})=(S^2A)',$$ which takes a cohomology
class to a cocycle from this class. This extension does not depend, up to an
isomorphism, on the choice of $\mu$ (compare Proposition 1.6) and possesses
the following partial (co-)universality property.\smallskip

{\sc Proposition} 2.6. {\it Let $M$ be an $A$-module with ${\frak m}M=0$.
Then the extension $(1)$ admits a unique homomorphism into an arbitrary 
extension $0\to M\to B\to A\to 0$ of $A$.}\smallskip

{\sc Proof.} The extension $0\to M\to B\to A\to0$ corresponds to some element 
of $H^2_{\rm Harr}(A;M)=H^2_{\rm Harr}(A;{\bb K})\otimes M$ (see Proposition 
2.2). The latter defines a mapping $H^2_{\rm Harr}(A;
{\bb K})'\to M$, which implies, in turn, a mapping $C\to B$. The resulting 
diagram 
$$
\def\mapright{\smash{\mathop{\longrightarrow}}}
\def\mapdown#1{\Big\downarrow\rlap{$\vcenter{\hbox{$\scriptstyle#1$}}$}}
\matrix{0&\mapright&H^2_{\rm Harr}(A;{\bb K})'&\mapright&C&\mapright&A&\mapright&0\cr
&&\mapdown{}&&\mapdown{}&&\mapdown{\rm id}&&\cr 0&\mapright&M&\mapright&B&
\mapright&A&\mapright&0\cr}$$
is an extension homomorphism. Its uniqueness is obvious. \medskip

{\bf 2.4.} $H^1_{\rm Harr}(A;M)$ is also interpreted as the set of automorphisms of 
any given extension $0\to M\, {\buildrel i\over\longrightarrow}\, B\, {\buildrel 
p\over\longrightarrow}\, A\to 0$ of $A$. An automorphism is an algebra automorphism
$f\colon B\to B$ such that $f\circ i=i$ and $p\circ f=p$. In previous notations
(see Proof of Proposition 2.3), $f(a,m)_q=(f_1(a,m),f_2(a,m))_q.$ The condition $p\circ
f=p$ means that $f_1(a,m)=a$. The condition $f\circ i=f$ means that $f_2(0,m)=m$, which
implies that $f_2(a,m)=f_2(a,0)+f_2(0,m)=m+\psi(a)$ (where $\psi(a)=f_2(a,0)$). The
multiplicativity of $f$ implies successively$$f((a_1,0)_q,(a_2,0)_q)=f(a_1,0)_qf(a_2,0)_q,$$
$$f(a_1a_2,\varphi_q(a_1,a_2))_q=(a_1,\psi(a_1))_q(a_2,\psi(a_2))_q,$$$$(a_1a_2,
\varphi_q(a_1,a_2)+\psi(a_1a_2))_q=(a_1a_2,\varphi_q(a_1,a_2)+a_1\psi(a_2)+a_2\psi(a_1)),$$$$
\psi(a_1a_2)=a_1\psi(a_2)+a_2\psi(a_1),$$ that is $d_1\psi=0$. Conversely, any $\psi\colon 
A\to M$ with $d_1\psi=0$ determines an algebra automorphism $f\colon B\to B, (a,m)_q\mapsto
(a,m+\psi(a))_q$ with the required properties.

Notice that $f(1,0)_q=(1,\psi(1))_q=(1,0)_q$, because $\psi(1)=0$ for
any derivation $\psi$.  Hence $f$ takes the unit element of $B$ into
the unit element of $B$ (cf. Proposition 2.4).  \medskip

{\bf 2.5.} In Section 4 we will use the following result due to
Harrison.\smallskip

{\sc Proposition 2.7.} ([Harr], Theorems 11 and 18). {\it Let $A={\bb K}[x_1,
\dots,x_n]$ be a polynomial algebra, and let $\frak m$ be the ideal of 
polynomials without constant terms. If an ideal $I$ of $A$ is contained in 
${\frak m}^2$, then
$$H^2_{\rm Harr}(A/I;{\bb K})\cong(I/({\frak m}\cdot I))'.$$}
Harrison's work contains also an explicit construction of the above 
homomorphism, which implies the following description of the canonical 
extension$$0\to H^2_{\rm Harr}(B;{\bb K})'\to C\to B\to0$$of $B=A/I$ (see 
2.3).\smallskip

{\sc Proposition 2.8.} {\it If $A,\ \frak m$, and $I$ are as in Proposition 
$2.7$, then the preceding extension for $B=A/I$ is$$0\to I/({\frak m}\cdot I)
\, {\buildrel i\over\longrightarrow}\, A/({\frak m}\cdot I)\, {\buildrel
p\over\longrightarrow}\, A/I\, \to0,$$where $i$ and $p$ are induced by the 
inclusions $I\to A$ and ${\frak m}\cdot I\to I$.}

\bigbreak

\centerline{\bf 3. Obstructions to extending deformations} \medskip

{\bf 3.1.} Let $\lambda$ be a deformation of a Lie algebra $L$ with a 
finite-dimensional local base $A$, and let $0\to {\bb K}\, {\buildrel i\over
\longrightarrow}\, B\, {\buildrel p\over\longrightarrow}\, A\to 0$ be a
1-dimensional extension of $A$, corresponding to a cohomology class $f\in
H^2_{\rm Harr}(A;{\bb K})$.

Let $I=i\otimes{\rm id}\colon L={\bb K}\otimes L\to B\otimes L$ and $P=
p\otimes{\rm id}\colon B\otimes L\to A\otimes L$. Let also $E=
\hat\varepsilon\otimes{\rm id}\colon B\otimes L\to {\bb K}\otimes L=L$, where
$\hat\varepsilon$ is the augmentation of $B$. The Lie algebra structure
$[\, ,\, ]_\lambda$ in $A\otimes L$ can be lifted to a $B$-bilinear 
operation $\{\, ,\, \}\colon\Lambda^2B\to B$ such that
\medskip

(i) $P\{l_1,l_2\}=[P(l_1),P(l_2)]_\lambda$ for any $l_1,l_2\in B\otimes
L$,

(ii) $\{I(l),l_1\}=I[l,E(l_1)]$ for any $l\in L, l_1\in B\otimes L$. 
\medskip
\noindent The operation $\{\, ,\, \}$ partially satisfies the Jacobi
identity, that is
$$\varphi(l_1,l_2,
l_3):=\{l_1,\{l_2,l_3\}\}+\{l_2,\{l_3,l_1\}\}+\{l_3,\{l_1,l_2\}\}\in{\rm
Ker}\, P.$$ Remark that $\varphi$ is multilinear and skew-symmetric, and
$\varphi (l_1,l_2,l_3)=0$ if $l_1\in{\rm Ker}\, E.$ (Indeed, if
$l_1=ml'_1$, where $m\in \hat{\frak m}={\rm Ker}\, \hat\varepsilon$,
then $\varphi(l_1,l_2,l_3)=\varphi
(ml'_1,l_2,l_3)=m\varphi(l'_1,l_2,l_3)=0.$) Hence $\varphi$ determines a
multilinear form$$\bar\varphi\colon\Lambda^3L=\Lambda^3((B\otimes
L)/{\rm Ker}\, E)\to {\rm Ker}\, P=L,$$ that is an element $\bar\varphi$
of $C^3(L;L)$. It is easy to check that $\delta\bar\varphi=0$.

Let $\{\, ,\, \}'$ be another $B$-bilinear operation $\Lambda^2B\to B$ satisfying 
the conditions (i), (ii) above. Then $\{l_1,l_2\}'-\{l_1,l_2\}\in{\rm Ker}\, P$
for any $l_1,l_2\in B\otimes L$, and if $l_1\in{\rm Ker}\, E$ then $\{l_1,l_2\}'
-\{l_1,l_2\}=0$ (as above, if $l_1=ml'_1,m\in\hat{\frak m}$, then $\{l_1,l_2\}'
-\{l_1,l_2\}=\{ml'_1,l_2\}'-\{ml'_1,l_2\}=m(\{l'_1,l_2\}'-\{l_1,l_2\})=0$.)
Hence the difference $\{\, ,\, \}'-\{\, ,\, \}$ determines a form 
$\psi\colon\Lambda^2L=\Lambda^2((B\otimes L)/{\rm Ker}\, E)\to{\rm Ker}\, P=L$,
that is determines a cochain $\psi\in C^2(L;L)$. Moreover, an arbitrary 
cochain $\psi\in C^2(L;L)$ may be obtained as $\{\, ,\, \}'-\{\, ,\, \}$ with
an appropriate $\{\, ,\, \}'$.

Using the cocycle $f_A$, it is easy to check that if $\bar\varphi,\bar\varphi'\in C^3(L;L)$ are the 
cochains corresponding to $\{\, ,\, \},\{\, , \}'$ in the sense of the 
construction above, then$$\bar\varphi'-\bar\varphi=\delta\psi.$$

Let ${\scr O}_\lambda(f)\in H^3(L;L)$ be the cohomology class of the cochain 
$\bar\varphi$. It is obvious that$${\scr O}_\lambda\colon H^2_{\rm Harr}(A;{\bb K})
\to H^3(L;L),\ f\mapsto {\scr O}_\lambda(f)$$is a linear map.

We can summarize the argumentation above in the following\smallskip

{\sc Proposition} 3.1. {\it The deformation $\lambda$ with base $A$ can
be extended to a  deformation of $L$ with base $B$ if and only if
${\scr O}_\lambda(f)=0$.} \smallskip

The cohomology class ${\scr O}_\lambda(f)$ is called the {\it obstruction} to
the extension of the deformation $\lambda$ from $A$ to $B$. \medskip

{\bf 3.2.} Suppose now that ${\scr O}_\lambda(f)=0$, that is the deformation 
$\lambda$ is extendible to a deformation with base $B$. We are going to 
study the set of all possible extensions.

Let $\mu,\mu'$ be deformations of $L$ with base $B$ such that $p_\ast\mu
=p_\ast\mu'=\lambda$. Then, according to 3.1, the difference $[\, ,\, ]_
{\mu'}-[\, ,\, ]_\mu$ determines and is determined by a certain cochain
$\psi\in C^2(L;L)$. Since $[\, ,\, ]_{\mu'}$ and $[\, ,\, ]_\mu$ both satisfy 
the Jacobi identity, $\delta\psi=0$. Moreover, it is easy to check that if we
replace any of the structures $[\, ,\, ]_\mu,[\, ,\, ]_{\mu'}$ with an 
equivalent one (see 1.1), then the cocycle $\psi$ will be replaced by a 
cohomologous cocycle. Thus the difference between two isomorphism classes of
deformations $\mu$ of $L$ with base $B$ such that $p_\ast\mu=\lambda$ is 
an arbitrary element of $H^2(L;L)$. In other words, $H^2(L;L)$ operates 
transitively on the set of these equivalence classes.

On the other hand, the group of automorphisms of the extension $0\to{\bb K}\,
{\buildrel i\over\longrightarrow}\, B\, {\buildrel p\over\longrightarrow}\,
A\to0$ also operates on the set of equivalence classes of deformations $\mu$. 
According to 2.5, this group is $H^1_{\rm Harr}(A;{\bb K})$, and according to 
Corollary 2.4, $H^1_{\rm Harr}(A;{\bb K})=({\frak m}/{\frak m}^2)'=TA$.
\smallskip

{\sc Proposition} 3.2. {\it These two operations are related to each other 
by the  
differential $d\lambda\colon TA\to H^2(L;L)$} (see Definition 1.10). {\it In
other words, if $r\colon B\to B$ determines an automorphism of the extension 
$0\to{\bb K}\, {\buildrel i\over\longrightarrow}\, B\, {\buildrel p\over
\longrightarrow}\, A\to0$ which corresponds to an element $h\in H^1_{\rm
Harr}(A;{\bb K})=TA$, then for any deformation $\mu$ of $L$ with base $B$ 
such that $p_\ast\mu=\lambda$, the difference between $[\, ,\, ]_{r_\ast\mu}$
and $[\, ,\, ]_\mu$ is a cocycle of the cohomology class $d\lambda(h)$.}
\smallskip

{\sc Proof} is obvious.\smallskip

{\sc Corollary} 3.3. {\it Suppose that the differential $d\lambda\colon
TA\to H^2(L;L)$ is onto. Then the group of automorphisms of the
extensions $0\to{\bb K}\, {\buildrel i\over\longrightarrow}\, B\,
{\buildrel p\over\longrightarrow}\, A\to0$ operates transitively on the
set of equivalence classes of deformations $\mu$ of $L$ with base $B$
such that $p_\ast\mu=\lambda$. In other words, $\mu$ is unique up to an
isomorphism and an automorphism of the extension $0\to{\bb K}\to B\to
A\to0$.} \medskip

{\bf 3.3.} The results of 3.1 and 3.2 may be generalized from the case 
of extension $0\to{\bb K}\to B\to A\to0$ to a more general case of extensions
$0\to M\, {\buildrel i\over\longrightarrow}\, B\, {\buildrel p\over
\longrightarrow}\, A\to0$, where $M$ is a finite-dimensional $A$-module
satisfying the condition ${\frak m}M=0$. The construction of 3.1 applied to 
a deformation $\lambda$ of $L$ with base $A$ yields an element of $H^3(L;
M\otimes L)=M\otimes H^3(L;L).$

The same element may be obtained from the previous construction in a
more direct way. Let $h\in M'$. We set $B_h=(B\oplus{\bb K})/{\rm Im}\,
(i\oplus h)$ (that is $B_h=B/i({\rm Ker}\, h)$ if $h\ne0$, and
$B_0=A\oplus{\bb K}$). There is an obvious extension $0\to{\bb K}\to
B_h\to A\to0$; let $f_h\in H^2_ {\rm Harr}(A;{\bb K})$ be the
corresponding cohomology class. The formula $h\mapsto{\scr
O}_\lambda(f_h)$ defines an element of ${\rm Hom}\,
(M',H^3(L;L))=M\otimes H^3(L;L)$ which coincides with the obstruction
constructed above.\smallskip

{\sc Proposition 3.4}. {\it A deformation $\mu$ of $L$ with base $B$
such that $p_\ast\mu=\lambda$ exists if and only if the element of $M\otimes
H^3(L;L)$ constructed above is equal to $0$. If $d\lambda\colon TA\to H^3
(L;L)$ is onto then the deformation $\mu$, if it exists, is unique up to an 
isomorphism and an automorphism of the extension $0\to M\to B\to A\to0$.}
\smallskip

{\sc Proof} is as above (see 3.1).
\bigbreak
\break			

\centerline{\bf 4. Construction of a miniversal deformation} \medskip

{\bf 4.1.} Suppose that $\dim H^2(L;L)<\infty.$
\smallskip

Let $C_0={\bb K},\, C_1={\bb K}\oplus H^2(L;L)'$, and let $$0\to
H^2(L;L)'\, {\buildrel i_1\over\longrightarrow}\, C_1\, {\buildrel
p'_1\over\longrightarrow}\, {\bb K}\to 0$$ be the canonical splitting
extension. The deformation $\eta_L$ of $L$ with base $C_1$ constructed
in 1.3 will be denoted here by $\eta_1$. Suppose that for some $k\ge1$
we have already constructed a finite-dimensional commutative algebra
$C_k$ and a deformation $\eta_k$ of $L$ with base $C_k$. Consider the
extension$$0\to H^2_{\rm Harr}(C_k;{\bb K})'\, {\buildrel \bar
i_{k+1}\over \longrightarrow}\, \bar C_{k+1}\, {\buildrel \bar
p'_{k+1}\over\longrightarrow} \, C_k\to0\eqno(2)$$ constructed in 2.3
using the cocycle $f_{C_k}$ (the notation was different there).
According to 3.3, we obtain the obstruction
$$
{\scr O}_{\eta_k}(f_{C_k}) \in H^2_{\rm Harr}(C_k,{\bb K})'\otimes H^3(L;L)
$$
to the extension of $\eta_k$. This gives us a map
$$\omega_k \colon H^2_{\rm Harr}(C_k,{\bb K})\to H^3(L;L).$$
Set
$$C_{k+1}=\bar C_{k+1}/\bar i_{k+1}\circ\omega_k'(H^3(L;L)').$$
Obviously, the extension (2) factorizes to an extension$$0\to({\rm
Ker}\, \omega_k)'\, {\buildrel i_{k+1}\over\longrightarrow}\, C_{k+1}\,
{\buildrel p'_{k+1}\over\longrightarrow}\, C_k\to 0.\eqno(3)$$

Notice that all the algebras $C_k$ are local.  Since $C_k$ is
finite-dimensional, the cohomology $H^2_{\rm Harr} (C_k;{\bb K})$ is
also finite-dimensional, and hence $C_{k+1}$ is
finite-dimensional.\smallskip

{\sc Proposition} 4.1. {\it The deformation $\eta_k$ admits an 
extension to a deformation with base $C_{k+1}$, and this extension is 
unique up to an isomorphism and an automorphism of an extension $(3)$.}
\smallskip

{\sc Proof}. According to Proposition 3.4, the obstruction to the
extension of the deformation $\eta_k$ of $L$ from $C_k$ to $C_{k+1}$ is
a homomorphism ${\rm Ker}\, \omega_k \to H^3(L;L)$, and it is easy to
show that it is precisely the restriction of $\omega_k$. Hence it is
equal to 0. The uniqueness of the extension is stated explicitly in
Proposition 3.4.

\medskip

We choose an extended deformation and denote it by $\eta_{k+1}$.

The induction yields a sequence of finite-dimensional algebras$${\bb K}
\, {\buildrel p'_1\over\longleftarrow}\, C_1\, {\buildrel p'_2\over
\longleftarrow}\, \dots\, {\buildrel p'_k\over\longleftarrow}\, C_k\, 
{\buildrel p'_{k+1}\over\longleftarrow}\, C_{k+1}\, {\buildrel p'_{k+2}
\over\longleftarrow}\, \dots,$$and a sequence of deformations $\eta_k$ of
$L$ such that $(p'_{k+1})_\ast\eta_{k+1}=\eta_k$.

Taking the projective 
limit, we obtain a formal deformation $\eta$ of $L$ with base 
$C=\overleftarrow{\displaystyle{\lim_{k\to\infty}}}C_k$.  In Theorem
4.5 below we will show that $\eta$ is a miniversal deformation of $L$. \medskip

{\bf 4.2.} Denote the space $H^2(L;L)$ briefly by $\bf H$.  Below we
assume that $\dim {\bf H} < \infty$.  Let ${\frak m}$ be the maximal
ideal in ${\bb K}[[{\bf H}']]$.

{\sc Proposition 4.2.} $C_k={\bb K}[[{\bf H}']]/I_k$ {\it where}
$${\frak m}^2=I_1\supset I_2\supset\dots,\ I_k\supset{\frak m}^{k+1}.$$
\smallskip
{\sc Proof.}  By construction,$$C_1={\bb K}\oplus{\bf H}'={\bb
K}[[{\bf H}']]/{\frak m}^2.$$ Suppose that we already know that
$$C_k={\bb K}[[{\bf H}']]/I_k,\ {\frak m}^2\supset I_k\supset{\frak m}^{k+1}.$$
Then, according to Proposition 2.8,
$$\bar C_{k+1}={\bb K}[[{\bf H}']]/({\frak m}\cdot I_k),$$ and by
construction $C_{k+1}$ is the quotient of $\bar C_{k+1}$ over an ideal
contained in $I_k/({\frak m}\cdot I_k)\subset{\frak m}^2/({\frak
m}\cdot I_k).$ Hence
$$C_{k+1}={\bb K}[[{\bf H}']]/I_{k+1},\ {\rm where}\ {\frak m}^2\supset
I_{k+1}\supset{\frak m}\cdot I_k\supset{\frak m}^{k+2}.$$
This completes the proof.\smallskip

{\sc Corollary 4.3.} {\it For $k\ge 2$ the projection $p_k' \colon C_k
\to C_{k-1}$ implies an isomorphism $TC_k \to TC_{k-1}$.  In
particular, the space $TC_k$ does not depend on $k$;
$TC_k=TC_1=H^2(L;L)$.  More precisely, for any $k \ge 1$ the
differential $d\eta_k\colon TC_k\to H^2(L;L)$ is an isomorphism.}
\smallskip

{\sc Proposition 4.4.} {\it $C={\bb K}[[{\bf H}']]/I$, where $I$ is an
ideal contained in ${\frak m}^2$.  Note that since ${\bb K}[[{\bf
H}']]$ is Noetherian, then $I$ is finitely generated.}\smallskip

{\sc Proof.} By construction, $C=\overleftarrow{\displaystyle{\lim_{k\to
\infty}}}C_k$ (see 4.1).  Proposition 4.2 gives an epimorphism
$$\overleftarrow{\lim_{k\to\infty}}{\bb K}[[{\bf H}']]/{\frak m}^{k+1}
	\to\overleftarrow{\lim_{k\to\infty}}C_k,$$
that is $${\bb K}[[{\bf H}']]\to C,$$
and 
$$C={\bb K}[[{\bf H}']]/I,\ {\rm  where}\ I=\cap I_k=\overleftarrow\lim\, I_k.
$$ \medskip

{\bf 4.3} {\sc Theorem 4.5.} {\it If $\dim H^2(L;L) < \infty$, then
the formal deformation $\eta$ is a miniversal formal deformation of
$L$.} 
\smallskip

{\sc Proof}.  Since $TC_k= H^2(L;L)$ and $d\eta_k={\rm id}$, then
$TC=H^2(L;L)$ and $d\eta={\rm id}$.  Let $A$ be a complete local
algebra with the maximal ideal $\frak m$, and let $\lambda$ be a
deformation of $L$ with base $A$.  We put $A_0=A/{\frak m}={\bb K}$
and $A_1=A/{\frak m}^2={\bb K}\oplus(TA)'$.  Then we fix a sequence of
1-dimensional extensions$$0\to{\bb K}\, {\buildrel
j_{k+1}\over\longrightarrow}\, A_{k+1}\, {\buildrel q_{k+1}\over
\longrightarrow}\, A_k\to0,\ k\ge1$$such that $A=\overleftarrow{\displaystyle
{\lim_{k\to\infty}}}A_k$. Let $Q_k\colon A\to A_k$ be the projection; we 
suppose that $Q_1$ is the natural projection $A\to A/{\frak m}^2$. Let 
$\lambda_k=(Q_k)_\ast\lambda$; it is a deformation of $L$ with base $A_k$. 
Obviously, $\lambda_k=(q_{k+1})_\ast\lambda_{k+1}$. We will construct 
inductively homomorphisms $\varphi_j\colon C_j\to A_j,\, j=1,2,\dots$ 
compatible with the projections $C_{j+1}\to C_j,\, A_{j+1}\to A_j$ and such 
that $(\varphi_j)_\ast\eta_j=\lambda_j$.

Define $\varphi_1\colon C_1\to A_1$ as ${\rm id}\oplus(d\lambda)'\colon{\bb K}
\oplus H^2(L;L)'\to{\bb K}\oplus(TA)'$; by definition of the differential, 
$(\varphi_1)_\ast\eta_1=\lambda_1$.
Suppose that $\varphi_k\colon C_k\to A_k$ with $(\varphi_k)_\ast\eta_k=
\lambda_k$ has been already constructed. The homomorphism $\varphi_k^\ast
\colon H^2_{\rm Harr}(A_k;{\bb K})\to H^2_{\rm Harr}(C_k;{\bb K})$ induced by
$\varphi_k$ takes the class of extension $0\to{\bb K}\to A_{k+1}\to A_k\to0$
into the class of some extension $0\to{\bb K}\to B\to C_k\to0$, and we have a
homomorphism
$$
\def\mapright{\smash{\mathop{\longrightarrow}}}
\def\mapdown#1{\Big\downarrow\rlap{$\vcenter{\hbox{$\scriptstyle#1$}}$}}
\matrix{0&\mapright&{\bb K}&\mapright&B&\mapright&C_k&\mapright&0\cr
&&\mapdown{}&&\mapdown{\psi}&&\mapdown{\varphi_k}&&\cr 
0&\mapright&{\bb K}&\mapright&A_{k+1}&\mapright&A_k&\mapright&0\cr}$$

Obviously, there exists a deformation $\xi$ of $L$ with base $B$ which 
extends $\eta _k$ (because the deformations $\lambda_k$ and $\eta_k$ have the 
same obstruction to extension) and such that $\psi_\ast\xi=\lambda_{k+1}$ 
(extensions of $\lambda_k$ and $\eta_k$ are both parameterized by $H^2(L;L)$).

According to Proposition 2.6, there exists a homomorphism

$$
\def\mapright#1{\smash{\mathop{\relbar\joinrel\longrightarrow}\limits^{#1}}}
\def\mapdown#1{\Big\downarrow\rlap{$\vcenter{\hbox{$\scriptstyle#1$}}$}}
\matrix{0&\mapright{}&H^2_{\rm Harr}(C_k;{\bb K})'&\mapright{\bar i_{k+1}}
&\bar C_{k+1}&\mapright{\bar p'_{k+1}}&C_k&\mapright{}&0\cr
&&\mapdown{r}&&\mapdown{\bar\chi}&&\mapdown{\rm id}&&\cr 
0&\mapright{}&{\bb K}&\mapright{}&B&\mapright{}&C_k&\mapright{}&0\cr}$$
and since the deformation $\eta_k$ is extended to $B$, it follows that the 
composition$$r \circ \omega_k'\colon H^3(L;L)'\to{\bb K}$$is zero. 
Hence the last diagram may be factorized to
$$
\def\mapright{\smash{\mathop{\longrightarrow}}}
\def\mapdown#1{\Big\downarrow\rlap{$\vcenter{\hbox{$\scriptstyle#1$}}$}}
\matrix{0&\mapright&({\rm Ker}\, \omega_k)'&\mapright&C_{k+1}&
\mapright&C_k&\mapright&0\cr
&&\mapdown{}&&\mapdown{\chi}&&\mapdown{\rm id}&&\cr 
0&\mapright&{\bb K}&\mapright&B&\mapright&C_k&\mapright&0\cr}$$
Since $d\eta_k\colon TC_k\to H^2(L;L)$ is an epimorphism (see 1.4.1), the two
deformations $\chi_\ast\eta_{k+1}$ and $\xi$ are related by some automorphism 
$f\colon B\to B$ of the extension $0\to{\bb K}\to B\to C_k\to0$. It remains to 
set $\varphi_{k+1}=\psi\circ f\circ\chi\colon C_{k+1}\to A_{k+1}$; indeed,
$(\varphi_{k+1})_\ast\eta_{k+1}=\psi_\ast\circ f_\ast\circ\chi_\ast\eta_{k+1}
=\psi_\ast\xi=\lambda_{k+1}$.

The limit map $\varphi\colon C\to A$ obviously satisfies the condition 
$\varphi_\ast\eta=\lambda$. The uniqueness property (ii) in Definition 
1.11 follows from the uniqueness in Proposition 1.8.  \medskip

{\bf 4.4} {\sc Theorem 4.6.} {\it If $\dim H^2(L;L) < \infty$, then the
base of the miniversal formal deformation of $L$ is formally embedded
in $H^2(L;L)$, that is, it may be described in $H^2(L;L)$ by a finite
system of formal equations.} \smallskip 

{\sc Proof.} Follows directly from Proposition 4.4. \smallskip

To make the computation of $C$ more specific, we need an appropriate
theory of Massey products.

\bigbreak

\centerline{\bf 5. Massey products} \medskip

{\bf 5.1.} The obstructions$$\omega_k \colon H^2_{\rm Harr}(C_k; {\bb
K})\to H^3(L;L)$$which arise in the construction of the miniversal
formal deformation of the Lie algebra $L$ (see 4.1) may be described in
terms of Massey products in $H^\ast(L;L)$. The appropriate theory of
Massey products was developed by the second author and Lang [FuL]. We
briefly recall this theory.\smallskip

{\sc Definition} 5.1. A differential graded Lie algebra (DGLA) is a vector 
space $\scr C$ over $\bb K$ with $\bb Z$ or ${\bb Z}_2$ grading ${\scr C}=
\bigoplus_i{\scr C}^i$ and with commutator operation $\mu\colon L\otimes L\to
L,\ \mu(\alpha\otimes\beta)=[\alpha,\beta]$ of degree 0 and a differential 
$\beta\colon{\scr C}\to{\scr C}$ of degree $+1$ satisfying the conditions
$$\eqalign{&[\alpha,\beta]=-(-1)^{\alpha\beta}[\beta,\alpha],\cr 
&\delta[\alpha,\beta]=[\delta \alpha,\beta]+(-1)^\alpha[\alpha,\delta \beta],
\cr &[[\alpha,\beta],\gamma]+(-1)^{\alpha(\beta+\gamma)}[[\beta,\gamma],
\alpha]+(-1)^{\gamma(\alpha+\beta)}[[\gamma,\alpha],\beta]=0,\cr}$$where the 
degree of a homogeneous element is denoted by the same letter as this element. 
\smallskip

Our main example of DGLA was introduced in 1.2: ${\scr C}^i=C^{i+1}(L;L)$.

The cohomology of $\scr C$ with respect to $\delta$ is denoted as ${\scr H}=
\bigoplus_i{\scr H}^i$. It is a graded Lie algebra. \medskip

{\bf 5.2.} The construction of Massey products in $\scr H$ given below
requires the following data. First, a graded cocommutative coassociative
coalgebra, that is a $\bb Z$ or ${\bb Z}_2$ graded vector space $F$ over 
$\bb K$ with a degree 0 mapping $\Delta\colon F\to F\otimes F$ 
(comultiplication) satisfying the conditions $S\circ\Delta=\Delta$, where 
$S\colon F\otimes F\to F\otimes F$ is defined as $S(\varphi\otimes\psi)=
(-1)^{\varphi\psi}(\psi\otimes\varphi)$, and $({\rm id}\otimes\Delta)\circ
\Delta=(\Delta\otimes{\rm id})\circ\Delta$. Second, a filtration $F_0\subset
F_1\subset F$ such that $F_0\subset{\rm Ker}\, \Delta$ and ${\rm Im}\, \Delta
\subset F_1\otimes F_1$.\smallskip

{\sc Proposition} 5.2 (see [FuL], Proposition 3.1). {\it Suppose a linear 
mapping $\alpha\colon F_1\to{\scr C}$ of degree 1 satisfies the condition
$$\delta\alpha=\mu\circ(\alpha\otimes\alpha) \circ\Delta.\eqno(4)$$Then
$$\mu\circ(\alpha\otimes\alpha)\circ\Delta(F) \subset{\rm Ker}\, \delta.$$}
 
(The right-hand side of the last formula is well defined because $\Delta(F)$ 
is contained in $F_1\otimes F_1$, the domain of $\alpha\otimes\alpha$).
\smallskip 
 
{\sc Definition} 5.3. Let $a\colon F_0\to{\scr H}$ and $b\colon F/F_1\to 
{\scr H}$ be linear maps of degrees 1 and 2. We say that $b$ is contained in the Massey $F$-product of 
$a$, and write $b\in\langle a\rangle_F$, or $b\in\langle a\rangle$, 
if there exists a degree 1 linear mapping $\alpha\colon F_1\to {\scr C}$ 
satisfying condition (4), and such that the diagrams 
 
$$ 
\def\mapright#1{\smash{\mathop{\relbar\joinrel\longrightarrow}\limits^{#1}}} 
\def\longmapright#1{\smash{\mathop{\relbar\joinrel\relbar\joinrel\relbar
\joinrel\relbar\joinrel\relbar\joinrel\rightarrow}\limits^{#1}}} 
\def\mapdown#1{\Big\downarrow\rlap{$\vcenter{\hbox{$\scriptstyle#1$}}$}} 
\matrix{F_0&\mapright{\alpha|_{F_0}}&{\rm Ker}\, \delta\cr 
\mapdown{\rm id}&&\mapdown{\pi}\cr F_0&\mapright{a}&{\scr H}\cr},\qquad\qquad 
\matrix{F&\longmapright{\mu\circ(\alpha\otimes\alpha)\circ\Delta}&{\rm Ker}\, 
\delta\cr \mapdown{\pi}&&\mapdown{\pi}\cr F/F_1&\longmapright{b}&{\scr H}\cr}
$$are commutative, where the vertical maps labeled by $\pi$ denote the 
projections of each space onto the quotient space. 
 
Note that the upper horizontal maps of the diagrams are well defined, 
since $\alpha(F_0)\subset\alpha({\rm Ker}\, \Delta)\subset{\rm Ker}\, \delta$ 
by virtue of (4), and $\mu\circ(\alpha\otimes\alpha)\circ\Delta(F) \subset{\rm 
Ker}\, \delta$ by Proposition 5.2. 

Note also that the definition makes sense even in the case, when 
$F_1=F$. In this case we do not need to specify any $b$, and we will simply 
say that $a$ {\it satisfies the condition of triviality of Massey 
$F$-products}. \smallskip 
 
{\sc Example} 5.4. Let $F$ be the dual of the maximal ideal of ${\bb K}[t]
/(t^{n+1})$, $F_0$ and $F_1$ be the duals of maximal ideals of ${\bb K}[t]/
(t^2)$ and ${\bb K}[t]/(t^n)$. Then $F_0$ and $F/F_1$ are 1-dimensional and 
are generated respectively by $t$ and $t^n$. In this case $a\colon F_0\to
\scr H$ and $b\colon F/F_1\to\scr H$ are characterized by $a(t)\in\scr H$ and
$b(t^n)\in\scr H$, and it is easy to check that $b\in\langle a\rangle_F$ if
and only if $b(t^n)$ belongs to the $n$-th Massey power of $a(t)$ in the
classical sense. In particular, for $n=2$, $b\in\langle a\rangle_F$ if and 
only if $b(t^2)=[a(t),a(t)]$. \medskip

{\bf 5.3.} The relationship between Massey products and Lie algebra 
deformations was established in the article [FuL] by the following result.

Let $A$ be a finite-dimensional local algebra with the maximal ideal $\frak
m$. Put $F=F_1={\frak m}'$ and $F_0=TA=({\frak m}/{\frak m}^2)'$.\smallskip

{\sc Proposition 5.5} ([FuL], Theorem 4.2). {\it A linear map $a\colon 
F_0\to H^2(L;L)$ is a differential of some deformation with base $A$ if and 
only if $-{1\over2}a$ satisfies the condition of triviality of Massey 
$F$-products}. \smallskip 

A similar result holds for formal deformations.

\bigbreak

\centerline{\bf 6. Calculating obstructions} \medskip

{\bf 6.1.} Adopt the notations of 4.1. Consider the sequence$${\bb K}\, 
{\buildrel p_1\over\longrightarrow}\, C_1\, {\buildrel p_2\over\longrightarrow}
\, \dots\, {\buildrel p_k\over\longrightarrow}\, C_k\, {\buildrel
\bar p_{k+1}\over\longrightarrow}\, \bar C_{k+1}.$$Recall that all $C_i,\bar 
C_i$ are finite-dimensional algebras,$$C_1={\bb K}\oplus H^2(L;L)',$$and
there is an extension$$0\to H^2_{\rm Harr}(C_k;{\bb K})'\, {\buildrel\bar 
i_{k+1}\over\longrightarrow}\, \bar C_{k+1}\, {\buildrel\bar p'_{k+1}\over
\longrightarrow}\, C_k\to0$$ and an obstruction homomorphism 
$$\omega_k\colon H^2_{\rm Harr}(C_k;{\bb K})\to H^3(L;L).$$Recall also that
$$C_{k+1}=\bar C_{k+1}/{\rm Im}\, (\bar i_{k+1}\circ\omega_k').$$

Let ${\frak m}_i,\bar{\frak m}_i$ be the maximal ideals in $C_i, \bar C_i$.  
Then we also have the sequence$${\frak m}_1\, {\buildrel p_2\over
\longrightarrow}\, {\frak m}_2\, {\buildrel p_3\over\longrightarrow}\, 
\dots\, {\buildrel p_k\over\longrightarrow}\,{\frak m}_k\, {\buildrel
\bar p_{k+1}\over\longrightarrow}\, \bar{\frak m}_{k+1}.$$ 
Consider the dual sequence$${\frak m}'_1\, {\buildrel p'_2\over
\longleftarrow}\, {\frak m}'_2\, {\buildrel p'_3\over\longleftarrow}\, 
\dots\, {\buildrel p'_k\over\longleftarrow}\,{\frak m}'_k\, {\buildrel
\bar p'_{k+1}\over\longleftarrow}\, \bar{\frak m}'_{k+1}.$$This is a sequence
of successively embedded cocommutative coassociative coalgebras. Put $\bar
{\frak m}'_{k+1}=F,\, {\frak m}'_1=F_0,\, {\frak m}'_k=F_1$. Then$$F_0=
H^2(L;L),\, F/F_1=H^2_{\rm Harr}(C_k;{\bb K}).$$We choose the grading in $F$ to
be trivial (${\rm deg}\, \varphi=0$ for any $\varphi\in F$). \medskip

{\bf 6.2.} {\sc Theorem 6.1.} $2\omega_k\in\langle
{\rm id}\rangle_F$ ({\it this inclusion refers to the Massey product in the 
sense of Definition $5.3$ in the cohomology ${\scr H}=\bigoplus_i{\scr H}^i,
\, {\scr H}^i=H^{i+1}(L;L)$, of the DGLA ${\scr C}=\bigoplus_i{\scr C}^i,\, 
{\scr C}^i=C^{i+1}(L;L))$. Moreover, an arbitrary element of $\langle{\rm 
id}\rangle_F$ is equal to $2\omega_k$ for an 
appropriate extension of the deformation $\eta_1=\eta_L$ of $L$ with 
base $C_1$ to a deformation $\eta_k$ of $L$ with base $C_k$.}\smallskip

{\sc Proof}. The Lie $C_k$-algebra structure $\eta_k$ on $C_k\otimes L$ is
determined by the commutators $[\l_1,l_2]_{\eta_k}\in C_k\otimes L$ of 
elements of $L=1\otimes L\subset C_k\otimes L$. The difference $[\, ,\, ]_
{\eta_k}-[\, ,\, ]$ is a linear map $\beta\colon\Lambda^2L\to{\frak m}_k
\otimes L$. This map may be regarded as a map ${\frak m}'_k=F_1\to{\rm Hom}
\, (\Lambda^2L,L)=C^2(L;L)$; we take the last map for $\alpha$ (see 
Definition 5.3). Obviously, $\alpha|F_0$ represents $a={\rm id}\colon F_0
\to H^2(L;L)$, and the Jacobi identity for $[\, ,\, ]_{\eta_k}$ means 
precisely that $\alpha$ satisfies condition (4). Moreover, it is clear, that 
different $\alpha$'s with these properties correspond precisely to different 
extensions $\eta_k$ of $\eta_1$. 

By definition, a map $b\colon F/F_1\to H^3(L;L)$ from $\langle a\rangle_F$ 
is represented by $\mu\circ(\alpha\otimes\alpha)\circ\Delta\colon F\to C^3
(L;L)$. On the other hand, the obstruction map $\omega_k\colon H^2_
{\rm Harr}(C_k;{\bb K})=\bar{\frak m}'_{k+1}/{\frak m}'_k=F/F_1\to H^3(L;L)$
is defined by means of lifting the commutator $[\, ,\, ]_{\eta_k}$ to a
skew-symmetric $\bar C_{k+1}$-bilinear operation $\{\, ,\, \}$ (satisfying
some additional conditions -- see 3.1). Choose a basis $m_1,\dots,m_s$ in
${\frak m}_k$, and extend it to a basis $\bar m_1,\dots,\bar m_s,\bar 
m_{s+1},\dots,$ $\bar m_{s+t}$ of $\bar{\frak m}_{k+1}$. (We will also consider 
the dual bases $\{m'_i\}$ and $\{\bar m'_i\}$ in ${\frak m}'_k$ and $\bar
{\frak m}'_{k+1}$.) Then$$[l_1,l_2]_{\eta_k}=[l_1,l_2]+\sum_{i=1}^sm_i
\otimes[l_1,l_2]_i,$$and the map $\alpha$ acts by the formula$$\alpha(m'_i)
(l_1,l_2)=[l_1,l_2]_i,\ i=1,\dots,,s.$$We define $\{\, ,\, \}$ by the formula
$$\{l_1,l_2\}=[l_1,l_2]+\sum_{i=1}^s\bar m_i\otimes[l_1,l_2]_i.$$Let the
multiplication in $\bar{\frak m}_{k+1}$ be$$\bar m_i\bar m_j=\sum_{p=1}
^{s+t}c_{ij}^p\bar m_p;$$then $\Delta\colon\bar{\frak m}'_{k+1}\to{\frak 
m}_k\otimes{\frak m}_k$ acts by the formula$$\Delta(\bar m'_p)=\sum_{i,j=1}
^sc_{ij}^pm'_i\otimes m'_j.$$We have$$\{\{l_1,l_2\},l_3\}=[l_1,l_2],l_3]+
\dots+\sum_{i,j=1}^s\sum_{p=s+1}^{s+t}c_{ij}^p\bar m_p\otimes[[l_1,l_2]_i,
l_3]_j,$$where ``$\dots$" denotes the part corresponding to $\bar m_1,\dots,
\bar m_s$. Thus the functional $\bar m'_p\in\bar{\frak m}'_{k+1}$ takes
$$\{\{l_1,l_2\},l_3\}+\{\{l_2,l_3\},l_1\}+\{\{l_3,l_1\},l_2\}$$into$$\sum
_{i,j=1}^sc_{ij}^p\left[\alpha(m'_j)(\alpha(m'_i)(l_1,l_2),l_3)+
\alpha(m'_j)(\alpha(m'_i)(l_2,l_3),l_1)+\alpha(m'_j)(\alpha(m'_i)(l_3,
l_1),l_2)\right]$$$$={1\over2}\mu\circ(\alpha\otimes\alpha)\circ\Delta(\bar 
m'_p),$$which shows that $\omega_k=\displaystyle{1\over2}b$. 
Theorem 6.1 follows.\bigbreak

\centerline{\bf 7. Further computations} \medskip

{\bf 7.1.} The goal of this Section is to provide a scheme of computation of 
the base of a miniversal deformation of a Lie algebra, convenient for 
practical use. We begin with the detailed description of the first two steps
of this inductive computation.

As in Section 4, we denote $H^2(L;L)$ by $\bf H$, and also denote by $\frak 
m$ the maximal ideal of the polynomial algebra ${\bb K}[{\bf H}']$. As 
before, we assume that $\dim{\bf H}<\infty$. Also we adopt the notations $C_k,
\bar C_k, {\frak m}_k, \bar{\frak m}_k$ of 4.1 and 6.1, and to avoid 
confusion, we denote the map $\alpha\colon{\frak m}'_k\to C^2(L;L)$ of 6.1 by 
$\alpha_k$.

According to 4.1,$$C_1={\bb K}\oplus{\bf H}'={\bb K}[{\bf H}']/{\frak m}^2,$$
and hence$${\frak m}_1={\frak m}/{\frak m}^2={\bf H}',\ {\frak m}'_1={\bf H}.
$$According to 4.2,$$\bar C_2={\bb K}[{\bf H}']/{\frak m}^3,$$and
hence$$ \bar{\frak m}_2={\frak m}/{\frak m}^3,\ \bar{\frak m}'_2 ={\bf
H}\oplus S^2{\bf H}.$$The map$$\alpha_1\colon{\frak m}'_1={\bf H}\to
C^2(L;L)$$takes a cohomology class into a representing cocycle. Hence
the map
$$\mu\circ (\alpha_1\otimes\alpha_1)\circ\Delta\colon\bar{\frak
m}'_2\to C^3(L;L), \eqno(5)$$ where $\Delta\colon\bar{\frak
m}'_2\to{\frak m}'_1\otimes{\frak m}' _1$ is the comultiplication,
acts as zero on $\bf H$ (because $\Delta|{\bf H} =0$) and takes
$\xi\eta\in S^2{\bf H}$ (where $\xi,\eta\in\bf H$) into the product of
the chosen cocycles $\alpha_1(\xi),\alpha_1(\eta)$ representing $\xi,
\eta$. Obviously (and according to Proposition 5.2), the image of the
map (5) belongs to ${\rm Ker}\, \delta$, and the composition of this
map with the projection $\pi\colon{\rm Ker}\, \delta\to H^3(L;L)$ acts
as zero on $\bf H$ and coincides with the multiplication $[\, ,\,
]\colon S^2{\bf H}\to H^3(L;L)$ on $S^2{\bf
H}$. Hence
$$\eqalign{{\frak m}'_2&={\bf H}\oplus{\rm Ker}\, ([\, ,\, ]\colon
	S^2{\bf H}\to H^3(L;L)),\cr 
{\frak m}_2&={{\frak m} \over{\frak m}^3+J_2},\ {\rm where}\ J_2={\rm
	Im}\, ([\, ,\, ]')\cr 
C_2&={{\bb K} [{\bf H}']\over{\frak m}^3+J_2}.\cr}$$
Note that if $\dim H^3(L;L)=q$, then $J_2$ is an ideal in ${\bb
K}[{\bf H}']$ generated by (at most) $q$ quadratic polynomials.

Furthermore, according to 4.2,$$\bar C_3={{\bb K}[{\bf H}']\over
{\frak m}^4+ ({\frak m}\cdot J_2)},$$and hence$$\bar{\frak
m}_3={{\frak m}\over {\frak m}^4+({\frak m}\cdot J_2)},\ \bar{\frak
m}'_3={\bf H}\oplus S^2{\bf H}\oplus K,$$where $K\subset S^3{\bf H}$
is the intersection of kernels of the maps
$$f_ \varphi\colon S^3{\bf H}\to H^3(L;L),\ \varphi\in{\bf H}',$$
$$f_\varphi(\xi\eta\zeta)=\varphi(\xi)[\eta,\zeta]+
	\varphi(\eta)[\zeta,\xi]+\varphi (\zeta)[\xi,\eta].$$
The map
$$\alpha_2\colon{\frak m}'_2={\bf H}\oplus{\rm Ker}\, [\, ,\, ]\to
C^2(L;L)$$ coincides with $\alpha_1$ on $\bf H$ and takes
$\sum\xi_i\eta_i\in{\rm Ker}\, [\, ,\, ]\ (\xi_i,\eta_i\in{\bf H})$
into a two-dimensional cochain whose coboundary is
$\sum[\alpha_1(\xi_i),\alpha_1 (\eta_i)]$. Hence the
composition
$$\mu\circ(\alpha_2\otimes\alpha_2)\circ
	\Delta\colon\bar{\frak m}'_3\to C^3(L;L),\eqno(6)$$
where $\Delta\colon\bar {\frak m}'_3\to{\frak m}'_2\otimes{\frak
m}'_2$ is the comultiplication, coincides with the map (5) on ${\bf
H}\oplus S^2{\bf H}$ and takes $\sum \xi_i\eta_i\zeta_i$
into$$[\alpha_1(\xi_i),\alpha_2(\eta_i,\zeta_i)]+
[\alpha_1(\eta_i),\alpha_2(\zeta_i,\xi_i)]+[\alpha_1(\zeta_i),\alpha_2
(\xi_i,\eta_i)].$$ According to Proposition 2.8, the latter is a
cocycle, and the composition of the map (6) and the projection
$\pi\colon{\rm Ker}\, \delta\to H^3(L;L)$ acts as zero on $\bf H$, as
$[\, ,\, ]$ on $S^2{\bf H}$, and as the ``triple Massey product" on
$K$. The kernel of this composition is 
${\frak m}_3'$, and ${\frak m}_3$ is the dual of this kernel.
Thus, by construction,
$${\frak m}_3={{\frak m}\over{\frak m}^4+J_3},\ C_3={{\bb K}[{\bf
	H}']\over{\frak m}^4+J_3},$$ 
where $J_3\cap S^2{\bf H}'= J_2\cap S^2{\bf H}'$. \medskip

{\bf 7.2.} Describe now the $k$-th induction step. Suppose that we have 
already constructed$$C_k={{\bb K}[{\bf H}']\over{\frak m}^{k+1}+J_k},\ {\frak
m}_k={{\frak m}\over{\frak m}^{k+1}+J_k},\ \alpha_k\colon{\frak m}'_k\to 
C^2(L;L).$$Then, according to 4.2,$$\bar C_{k+1}={{\bb K}[{\bf H}']\over
{\frak m}^{k+2}+({\frak m}\cdot J_k)},\ \bar{\frak m}_{k+1}={{\frak m}\over
{\frak m}^{k+2}+({\frak m}\cdot J_k)},$$$$\bar{\frak m}'_{k+1}\subset{\bf H}
\oplus S^2{\bf H}\oplus\dots\oplus S^{k+1}{\bf H}.$$The image of the 
composition$$\mu\circ(\alpha_k\otimes\alpha_k)\circ\Delta\colon\bar{\frak m}'
_{k+1}\to C^3(L;L),\eqno(7)$$where $\Delta\colon\bar{\frak m}'_{k+1}\to
{\frak m}'_k\otimes{\frak m}'_k$ is the comultiplication, is contained in 
${\rm Ker}\, \delta$ (Proposition 5.2), and the composition$$\pi\circ\mu\circ 
(\alpha_k\otimes\alpha_k)\circ\Delta\colon\bar{\frak m}'_{k+1}\to 
H^3(L;L)$$acts as zero on ${\frak m}'_k$. We put$${\frak m}'_{k+1}={\rm Ker}
\, (\pi\circ\mu\circ(\alpha_k\otimes\alpha_k)\circ\Delta)\supset{\frak m}'
_k.$$The map $\alpha_k\colon{\frak m}'_k\to C^2(L;L)$ is extended to the map 
$\alpha_{k+1}\colon{\frak m}'_{k+1}\to C^2(L;L)$ such that $\delta\circ
\alpha_{k+1}$ is the restriction of the map (7). The dual to ${\frak m}'_
{k+1}$ is$${\frak m}_{k+1}={{\frak m}\over{\frak m}^{k+2}+J_{k+1}},$$and we 
put$$C_{k+1}={\bb K}\oplus{\frak m}_{k+1}={{\bb K}[{\bf H}']\over{\frak m}^
{k+2}+J_{k+1}}.$$This completes the construction. \medskip

{\bf 7.3.} Two following useful observations are easily derived from the 
description of the construction given in 7.1 -- 7.2.\smallskip

{\sc Proposition 7.1}. {\it For} $l\le k$,$$J_{k+1}\cap S^l{\bf H}'=J_k\cap S^l
{\bf H}'.$$

{\sc Proof}. We use induction with respect to $k$. For $k=2$ this was 
proved in 7.1. Suppose that $J_k\cap S^{l-1}{\bf H}'=J_{k-1}\cap S^{l-1}{\bf
H}'$. Then $({\frak m}\cdot J_k)\cap S^l{\bf H}'=({\frak m}\cdot J_{k+1})\cap
S^l{\bf H}'$. Hence $\bar{\frak m}'_{k+1}$ and $\bar{\frak m}'_k$ have the 
same $S^l{\bf H}'$ component. Since $\Delta$ has degree 0 with respect to 
${\bf H}'$, and $\alpha_k$ coincides with $\alpha_{k-1}$ on ${\frak m}'_{k-1}
$, we may conclude that ${\frak m}'_{k+1}$ and ${\frak m}'_k$ also have the 
same $S^l{\bf H}'$ component. Proposition 7.1 follows.\smallskip

{\sc Proposition 7.2.} {\it If $\dim H^3(L;L)=q$, then the ideal $I=
\overleftarrow{\lim}\, I_k=\overleftarrow{\lim}\, J_k$ from 
Proposition~$4.4$ has at most $q$ generators. Less formally, the base
of the miniversal deformation of $L$ is the zero locus of a formal map
$H^2(L;L)\to H^3(L;L)$.}\smallskip

{\sc Proof.} By construction, ${\frak m}_k=\bar{\frak m}_k/G_k$, where $G_k$
is generated by the image of a certain map $\beta_k\colon H^3(L;L)'\to\bar
{\frak m}_k$ (namely, $\beta_k=(\pi\circ\mu\circ(\alpha_{k-1}\otimes\alpha
_{k-1})\circ\Delta)'$, see 7.2). Actually, $\bar{\frak m}_k$ is a quotient of
$\bar{\frak m}_{k+1}$, and $\beta_{k+1}$ is a lift of $\beta_k$. Put$$\bar
{\frak m}_\infty=\overleftarrow{\lim_{k\to\infty}}\bar{\frak m}_k,\  
{\frak m}_\infty=\overleftarrow{\lim_{k\to\infty}}{\frak m}_k.$$Then ${\frak
m}_\infty=\bar{\frak m}_\infty/G_\infty$, where $G_\infty$ is generated by 
the image of$$\beta_\infty=\overleftarrow{\lim_{k\to\infty}}\beta_k\colon
H^3(L;L)'\to\bar{\frak m}_\infty.$$Furthermore,$${\frak m}_\infty={\frak m}/I,
\ \bar{\frak m}_\infty={\frak m}/({\frak m}\cdot I),$$where $I=\overleftarrow
{\lim}\, I_k$. Hence$$G_\infty=I/({\frak m}\cdot I),$$and it is clear that
generators of $G_\infty$ are lifted to generators of $I$. Since $G_\infty$
is generated by (at most) $q$ generators, Proposition 7.2 follows. \medskip

{\bf 7.4.} We conclude this section with a brief discussion of the graded 
case. Suppose that the Lie algebra $L$ is $G$-graded, where $G$ is an Abelian 
group: $L=\bigoplus_{g\in G}L_q,\ [L_g,L_h]\subset L_{g+h}.$ In this case the
cochains $\scr C$ and the cohomology $\scr H$ get an additional grading:
$C^q(L;L)=\bigoplus_{g\in G}C^q_{(g)}(L;L)$ $(\varphi\in C^q_{(g)}(L;L)$, if
$\varphi(l_1,\dots,l_q)\in L_{g_1+\dots+g_q-g}$ for $l_1\in L_{g_1},\dots,
l_q\in L_{g_q})$, and $H^q(L;L)=\bigoplus_{g\in G}H^q_{(g)}(L;L)$. The 
condition $\dim H^2(L;L)<\infty$ may be replaced in this case by a weaker 
condition: $\dim H^2_{(g)}(L;L)<\infty$ for each $g$. We preserve the notation 
$\bf H$ for $H^2(L;L)$, but ${\bf H}'$ will denote $\bigoplus_{g\in G}H^2_{(g)}
(L;L)'$. All the spaces ${\bf H}, {\bf H}', C_k, {\frak m}_k, \bar{\frak m}_k,
{\frak m}'_k, \bar{\frak m}'_k$ have natural $G$-gradings, and all the maps 
$\alpha_k$ have degree 0. The whole construction is modified correspondingly. 
We restrict ourselves to the modified version of Proposition 7.2.\smallskip

{\sc Proposition 7.3.} {\it The ideal $I$ from Proposition $4.4$ is always 
generated by homogeneous elements. Moreover, if $\dim H^3_{(g)}(L;L)=q_g$, 
then $I$ has at most $q_g$ generators of degree $g$. Less formally, the base
of the miniversal deformation of $L$ is the zero locus of a formal map $H^2
(L;L)\to H^3(L;L)$ of degree $0$.}\bigbreak

\centerline{\bf 8. Example: Deformations of the Lie algebra ${\bi L}_{\bf 
1}$} \medskip

{\bf 8.1.} Let $L_1$ be the complex Lie algebra of polynomial vector 
fields $p(x)\displaystyle{d\over dx}$ on the line such that 
$p(0)=p'(0)=0$. The deformations of this Lie algebra were studied by 
the first author ([Fi2], [Fi3]), and its formal miniversal deformation 
was completely described in our joint paper [FiFu]. It turned out that 
geometrically the base of this deformation is the union of three 
algebraic curves with a common point: two non-singular, having a 
common tangent, and one with a cusp, where the tangent at the cusp 
coincides with the tangent to the smooth components. 

Below we show how these results can be obtained by the methods of this 
article. We will need some (surprisingly little) information about the 
cohomology and deformations of the Lie algebra $L_1$. All this 
information is contained in the articles [FeFu], [Fi2], [Fi3], 
[FiFu]. \medskip

{\bf 8.2.} As a complex vector space, the Lie algebra $L_1$ has the basis 
$\{e_i|i\ge1\},\ e_i=x^{i+1}\displaystyle{d\over dx}$, 
and the commutator operation is 
$[e_i,e_j]=(j-i)e_{i+j}$. This Lie algebra is $\bb Z$-graded, ${\rm deg}\,
e_i=i$.\smallskip

{\sc Proposition 8.1} ([FeFu], [Fi2]). {\it The dimensions of $H^2(L_1;L_1)$ 
and $H^3(L_1;L_1)$ are equal to $3$ and $5$. Moreover,$$\eqalign{\dim H^2_
{(q)}(L_1;L_1)&=\cases{1&if $2\le q\le4$,\cr 0&otherwise;\cr}\cr \dim H^3_
{(q)}(L_1;L_1)&=\cases{1&if $7\le q\le11$,\cr 0&otherwise.\cr}\cr}$$}
\smallskip

{\sc Proposition 8.2} ([Fi3]). {\it Let $0\ne\alpha\in H^2_{(2)}(L_1;L_1), 
0\ne\beta\in H^2_{(3)}(L_1;L_1), 0\ne\gamma\in H^2_{(4)}(L_1;L_1).$ Then 
$0\ne[\beta,\gamma]\in H^3_{(7)}(L_1;L_1), 0\ne[\gamma,\gamma]\in H^3_{(8)}
(L_1;L_1)$. Furthermore, $0\ne\langle\beta,\beta,\beta\rangle\in H^3_{(9)}
(L_1;L_1).$}\smallskip

The latter means that if $b\in C^2_{(3)}(L_1;L_1)$ is a representative of 
$\beta$, and if $[b,b]=\delta g,\, g\in C^2_{(6)}(L_1;L_1)$, then the 
cohomology class of the cocycle $[b,g]\in C^3_{(9)}(L_1;L_1)$ (which does not 
depend on the choice of $b$ and $g$) is not equal to 0. \medskip

{\bf 8.3.} Here are some explicit constructions of deformations of the
Lie algebra $L_1$.\smallskip
 
{\sc Proposition 8.3} ([Fi2]). {\it The formulas
$$\eqalign{[e_i,e_j]^1_t&=(j-i)(e_{i+j}+te_{i+j-1});\cr
[e_i,e_j]^2_t&=\cases{(j-i)e_{i+j}&if $i\ne1,j\ne1$,\cr
(j-1)e_{j+1}+tje_j&if $i=1,j\ne1;$\cr}\cr
[e_i,e_j]^3_t&=\cases{(j-i)e_{i+j}&if $i\ne2,j\ne2$,\cr
(j-2)e_{j+2}+tje_j&if $i=2,j\ne2$\cr}\cr}$$ determine three
one-parameter deformations of the Lie algebra $L_1$. All the three
deformations are pairwise not equivalent. Moreover, if
$L_1^1,L_1^2,L_1^3$ are Lie algebras from the three families
corresponding to arbitrary non-zero values of the parameter (up to an
isomorphism, they do not depend on the non-zero parameter value), then
neither two of $L_1^1,L_1^2,L_1^3$ are isomorphic to each other.}
\smallskip

{\sc Corollary 8.4.} {\it The base of any versal deformation of the Lie 
algebra $L_1$ contains at least three different irreducible curves.} \medskip

{\bf 8.4.} We will use the notations of Section 7. Let $\alpha,\beta,\gamma$ 
be a basis of ${\bf H}=H^2(L_1;L_1)$ (as in Proposition 8.2), and let $x,y,z$ 
be the dual basis in ${\bf H}'$. The algebra $S^\ast{\bf H}'={\bb C}[x,y,z]$ 
has the monomial basis $\{x^py^qz^r\}$. Let $\{\alpha^p\beta^q\gamma^r\}$ be 
the dual basis in the coalgebra $S^\ast{\bf H}$; the comultiplication 
$\Delta\colon S^\ast{\bf H}\to S^\ast{\bf H}\otimes S^\ast{\bf H}$ acts by 
the formula$$\Delta(\alpha^p\beta^q\gamma^r)=\sum_{i=0}^p\sum_{j=0}^q\sum_
{k=0}^r\alpha^i\beta^j\gamma^k\otimes \alpha^{p-i}\beta^{q-j}\gamma^{r-k}.$$

Choose cocycles $a\in C^2_{(2)}(L_1;L_1),\ b\in C^2_{(3)}(L_1;L_1),\ c\in 
C^2_{(4)}(L_1;L_1)$ representing $\alpha,\beta,\gamma$. Then$$\alpha_1\colon
{\frak m}'_1={\bf H}\to C^2(L_1;L_1)$$is defined by the formulas$$\alpha_1
(\alpha)=a,\ \alpha_1(\beta)=b,\ \alpha_1(\gamma_)=c.$$
Since $H^3_{(q)}(L_1;L_1)=
0$ for $q<7$, there exist $d\in C^2_{(4)}(L_1;L_1), e\in C^2_{(5)}
(L_1;L_1), f\in C^2_{(6)}(L_1;L_1), g\in C^2_{(6)}(L_1;L_1)$, such that 
$[a,a]=\delta d,[a,b]=\delta e,[a,c]=\delta f,[b,b]=\delta g$ (the notation 
$g$ has been already used in 8.2). Since $c\in C^2_{(4)}(L_1;L_1)$ is a 
cocycle, we can replace $d$ with $d+tc$, where $t$ is an arbitrary complex
number. Finally, since $\delta[a,d]=0$, we also have $[a,d]=\delta h$ for 
some $h\in C^2_{(6)}(L_1;L_1).$

The space $\bar{\frak m}'_2={\bf H}\oplus S^2{\bf H}$ is spanned by $\alpha, 
\beta,\gamma,\alpha^2,\alpha\beta,\alpha\gamma,\beta^2,\beta\gamma,\gamma^2$.
The map $\mu\circ(\alpha_1\otimes\alpha_1)\circ\Delta\colon\bar{\frak m}'_2
\to C^3(L_1;L_1)$ acts in the following way:$$\alpha,\beta,\gamma\mapsto0;\ 
\alpha^2\mapsto\delta d,\ \alpha\beta\mapsto2\delta e,\ 
\alpha\gamma\mapsto2\delta
f,\ \beta^2\mapsto\delta g,$$$$\beta\gamma\mapsto2[b,c]\notin{\rm Im}\, \delta,
\gamma^2\mapsto[c,c]\notin{\rm Im}\, \delta.$$Hence ${\frak m}'_2$ is 
generated by $\alpha,\beta,\gamma,\alpha^2,\alpha\beta,\alpha\gamma,\beta^2,$
and$$\alpha_2\colon{\frak m}'_2\to C^2(L_1;L_1)$$ is defined as
$\alpha_1$ on ${\bf H}$ and 
$$\alpha_2(\alpha^2)=d+tc,\ \alpha_2(\alpha\beta)=2e,\ 
\alpha_2(\alpha\gamma)=2f,\ \alpha_2(\beta^2)=g.$$
Furthermore,$${\frak m}_2
={{\frak m}\over{\frak m}^3+(yz,z^2)},$$$$\bar{\frak m}_3={{\frak m}\over{\frak m}
^4+({\frak m}\cdot(yz,z^2))}={{\frak m}\over{\frak m}^4+(xyz,xz^2,y^2z,yz^2,
z^3)},$$and$$\bar{\frak m}'_3={\bf H}\oplus S^2{\bf H}\oplus K,$$where $K$ 
is the subspace of $S^3{\bf H}$ spanned by $\alpha^3,\alpha^2\beta,\alpha^2
\gamma,\alpha\beta^2,\beta^3$. The map $\mu\circ(\alpha_2\otimes\alpha_2)\circ
\Delta\colon\bar{\frak m}'_3\to C^3(L_1;L_1)$ acts as $\mu\circ(\alpha_1
\otimes\alpha_1)\circ\Delta$ on ${\bf H}\oplus S^2{\bf H}$ (see above), and 
acts on $K$ in the following way:$$\eqalign{\alpha^3&\mapsto2[a,d+tc]=\delta
(h+tf),\cr \alpha^2\beta&\mapsto4[a,e]+2[b,d]+2t[b,c],\cr \alpha^2\gamma
&\mapsto4[a,f]+2[c,d]+2t[c,c],\cr \alpha\beta^2&\mapsto2[a,g]+4[b,e],\cr 
\beta^3&\mapsto2[b,g]\notin{\rm Im}\, \delta.\cr}$$Since $4[a,e]+2[b,d]\in
{\rm Ker}\, \delta,\ [b,c]\notin{\rm Im}\, \delta,$ and $\dim H^3_{(7)}(L_1;
L_1)=1,$ we can choose $t$ in such a way that the image of $\alpha^2\beta$ is 
cohomologous to 0, $$\alpha^2\beta\mapsto\delta k,\ k\in C^2_{(7)}(L_1;L_1).$$
Since $4[a,f]+2[c,d]+2t[c,c],\, 2[a,g]+4[b,e]\in{\rm Ker}\, \delta,\ [c,c]
\notin{\rm Im}\, \delta$, and $\dim H^3_{(8)}(L_1;L_1)=1$, there exist 
complex numbers $A,B$ such that the images of $\alpha^2\gamma-A\gamma^2,\
\alpha\beta^2-B\gamma^2$ are cohomologous to 0,$$\eqalign{\alpha^2\gamma-
A\gamma^2&\mapsto\delta l,\ l\in C^2_{(8)}(L_1;L_1),\cr \alpha\beta^2-
B\gamma^2&\mapsto\delta m,\ m\in C^2_{(8)}(L_1;L_1).\cr}$$Hence ${\frak 
m}'_3$ is generated by the generators of ${\frak m}'_2$ (see above) and
also $\alpha^3,\alpha^2\beta,\alpha^2\gamma-A\gamma^2,$ $\alpha\beta^2-
B\gamma^2$. Thus$${\frak m}_3={{\frak m}\over{\frak m}^4+(yz,z^2+Ax^2z+
Bxy^2,y^3)}.$$

To complete this description of the base of the miniversal deformation of 
$L_1$, we need to continue the induction to calculate ${\frak m}_4$ and
${\frak m}_5$. This would require more information about the multiplications 
in the cohomology of $L_1$. It turns out, however, that we can avoid any 
additional computations if we use Corollary 8.4. \medskip

{\bf 8.5.} According to Propositions 7.3 and 8.1, the base of the miniversal 
deformation of $L_1$ is ${\bb C}[[x,y,z]]/(F_1,F_2,F_3,F_4,F_5)$, where $F_1,
\dots,F_5$ are polynomials in $x,y,z$ of degrees $7,\dots,11$ (with ${\rm deg}
\, x=2, {\rm deg}\, y=3, {\rm deg}\, z=4$). The calculations of 8.4 show that
$$\eqalign{F_1&=yz+\dots,\cr F_2&=z^2+Ax^2z+Bxy^2+\dots,\cr F_3&=y^3+\dots,
\cr}$$where ``$\dots$'', and $F_4,F_5$ as well, are linear combinations of 
4- and 5-fold products of $x,y,z$ having appropriate degrees. These products
are the following monomials.$$\eqalign{{\rm degree}\ 7&\colon\ {\rm none},\cr
{\rm degree}\ 8&\colon\ x^4, \cr {\rm degree}\ 9&\colon\ x^3y,\cr {\rm 
degree}\ 10&\colon\ x^5,x^3z,x^2y^2, \cr {\rm degree}\ 11&\colon\ x^4y,x^2yz.
\cr}$$We exclude the monomial $x^2yz$, because it can be extinguished by 
adding a constant times $x^2F_1$, and get the following intermediate result.
\smallskip

{\sc Lemma 8.5.} {\it The base of the miniversal deformation of $L_1$ is 
described in $H^2(L_1;L_1)$ by a system of formal equations$$\eqalign{&
\beta\gamma=0,\cr &\gamma^2+A\alpha^2\gamma+B\alpha\beta^2+C\alpha^4=0,\cr
&\beta^3+D\alpha^3\beta=0,\cr &E\alpha^5+F\alpha^3\gamma+G\alpha^2\beta^2=0,
\cr &H\alpha^4\beta=0.\cr}\eqno(8)$$}\smallskip

Consider the zero locus $X$ of the first three equations (8).\smallskip

{\sc Lemma 8.6}. {\it If $C=BD,\, A^2\ne 4C,$ and $D\ne0$, then $X$ is the 
union of three irreducible curves. Otherwise $X$ does not contain three 
different irreducible curves.}\smallskip

{\sc Proof.} Let $(\alpha,\beta,\gamma)\in X$. The first equation (8) says 
that either $\beta=0$, or $\gamma=0$. If $\beta=0$, then the third equation
holds, and the second equation becomes$$\gamma^2+A\alpha^2\gamma+C\alpha^4=
(\gamma+u\alpha^2)(\gamma+v\alpha^2)=0,\eqno(9)$$where $u\ne v$ if $A^2\ne
4C$. Hence $X\cap\{\beta=0\}$ is the union of two parabolas. If $\gamma=0$,
then the second and the third equations become$$\eqalign{\alpha(B\beta^2+
C\alpha^3)&=0,\cr \beta(\beta^2+D\alpha^3)&=0,\cr}$$which describes just one 
point $\alpha=0,\beta=0$ if $C\ne BD$, the semicubic parabola $\beta^2+D
\alpha^3=0$ if $0\ne C=BD$, and the union of the same semicubic parabola and 
the line $\beta=0$ if $0=C=BD$. In the last case one of the curves (9) is also
the line $\beta=0,\gamma=0$. Lemma 8.6 follows.\smallskip

{\sc Theorem 8.7}. {\it The base of the miniversal deformation of the Lie 
algebra $L_1$ is described in $H^2(L_1;L_1)$  by the system of formal 
equations$$\eqalign{&\beta\gamma=0,\cr &\gamma^2+A\alpha^2\gamma+B\alpha
(\beta^2+D\alpha^3)=0,\cr &\beta(\beta^2+D\alpha^3)=0,\cr}$$where $A^2\ne4BD,$ 
and $D\ne0$.}\smallskip

{\sc Proof.} Corollary 8.4 and Lemma 8.6 imply that in equations (8) $C=BD,
A^2\ne4C$, and $D\ne0$. Hence the three curves, which are contained in the 
base of the miniversal deformation according to Corollary 8.4, are$$\eqalign{
\beta=0&,\ \gamma+u\alpha^2=0;\cr \beta=0&,\ \gamma+v\alpha^2=0;\cr \gamma=0
&,\ \beta^2+D\alpha^3=0,\cr}$$where $u\ne v,\, u+v=A,\, uv=BD$. Hence the 
left hand sides of the last two equations (8) should be equal to 0 on these 
curves. The monomial $\alpha^4\beta$ is not equal to 0 on the third of the 
curves; hence $H=0$. If $\beta=0$, then the fourth equation becomes $\alpha^3
(E\alpha^2+F\gamma)=0$, which cannot hold on {\it both} parabolas $\gamma+
u\alpha^2=0,\, \gamma+v\alpha^2=0$ unless $E=F=0$. Finally, if $\gamma=0$, 
then the fourth equation (with $E=F=0$) becomes $G\alpha^2\beta^2=0$ which 
does not hold on the third curve unless $G=0$. \medskip

{\bf 8.6.} Note that the computations made in the article [FiFu] let us find 
the constants $A,B,D$ from Theorem 8.7. Since these constants depend on a 
particular choice of cocycles $a\in C^2_{(2)}(L_1;L_1),\ b\in C^2_{(3)}(L_1;
L_1),\ c\in C^2_{(4)}(L_1;L_1)$ representing generators of $H^2(L_1;L_1)$, we 
need to specify these cocycles first.

Let $W$ be an $L_1$-module spanned by $e_j$ with all $j\in\bb Z$ and with the 
$L_1$-action $e_i(e_j)=(j-i)e_{i+j}$. It is an extension of the adjoint 
representation. Define a cochain $$\mu_k\in C^1_{(k)}(L_1;W),\ k\ge2,$$by the 
formula$$\mu_k(e_i)=(-1)^{i+1}{k-1\choose i-2}e_{i-k}.$$\smallskip

{\sc Proposition 8.8} [FiFu]. {\it If $k=2,3,4$, then $\delta\mu_k$ belongs 
to $C^2_{(k)}(L_1;L_1)$ and is a cocycle not cohomologous to $0$.}\smallskip

{\sc Proposition 8.9} [FiFu]. {\it If one chooses $a,b,c$ to be $\delta\mu_2,
\, \delta\mu_3,\, \delta\mu_4$, then}$$A=-{2\cdot11\cdot37\over5\cdot13^2},\ 
B={4\cdot7\cdot17\over3\cdot25\cdot13},\ D={32\cdot27\over13^3}.$$

\vskip 20pt

\centerline{\bf BIBLIOGRAPHY}

\bigskip
\frenchspacing
\parindent=50pt

\bib{B}Barr, M., ``Harrison homology, Hochschild homology and 
triples,'' {\it J. Algebra}, {\bf 8} (1963), 314--323.

\medskip

\bib{FeFu}Feigin, B., Fuchs, D., ``Homology of the Lie algebras of 
vector fields on the line,'' {\it Funct. Anal. Appl.}, {\bf 14:3} (1980),
45--60.

\medskip

\bib{Fi1}Fialowski, A., ``Deformations of Lie Algebras,'' {\it Math. USSR-
Sbornik,} {\bf 55:2} (1986), 467--473.

\medskip

\bib{Fi2}Fialowski, A., ``Deformations of the Lie algebra of vector 
fields on the line,'' {\it Russian Math. Surveys}, {\bf 38} (1983), 185--186.

\medskip

\bib{Fi3}Fialowski, A., ``An example of formal deformations of Lie 
algebras,'' {\it NATO Conference on Deformation Theory of Algebras and 
Applications}, Il Ciocco, Italy, 1986. Proceedings, Kluwer, 
Dordrecht, 1988, 375--401. 

\medskip

\bib{FiFu}Fialowski, A., Fuchs, D., ``Singular deformations of Lie 
algebras on an example,'' in {\it Topics in Singularity Theory:
V.~I.~Arnold's {\rm 60}${}^{th}$ Anniversary Collection,} A.~Khovanskii,
A.~Varchenko, V.~Vassiliev (Eds.), Transl.~A.M.S.~Ser.~2,
{\bf 180}, pp.~77--92, Amer.~Math.~Soc., Providence, RI, 1997.

\medskip

\bib{Fu}Fuchs, D., {\it Cohomology of infinite-dimensional Lie 
Algebras}, Consultants Bureau, NY, London (1986).

\medskip

\bib{FuL}Fuchs, D., Lang, L., ``Massey products and Deformations,'' 
To appear in {\it Journal of Pure and Applied Algebra}.

\medskip

\bib{G}Gerstenhaber, M., ``On the deformation of rings and algebras,
I, III,'' {\it Ann. Math.} {\bf 79} (1964), 59--103; {\bf 88} (1968), 1--34.

\medskip

\bib{GoM}Goldman, W.~M., Millson,~J.~J., ``The Deformation Theory of
Representations of Fundamental Groups of Compact K\"ahler
Manifolds,'' {\it IHES Pub. Math.}, {\bf 67} (1988), 43--96.

\medskip

\bib{Harr}Harrison, D.K., ``Commutative algebras and cohomology,'' 
{\it Trans. Amer. Math. Soc.}, {\bf 104} (1962), 191--204. 

\medskip

\bib{Hart}Hartshorne, R., {\it Algebraic Geometry}, Springer 
(1977).

\medskip

\bib{I}Illusie, L., ``Complexe cotangent et d\'eformations I,'' {\it 
Lect Notes in Math.} {\bf 239}, Springer (1971). 

\medskip

\bib{K}Kontsevich, M., {\it Topics in Algebra: Deformation Theory}, 
Lecture Notes, Univ. Calif. Berkeley (1994).

\medskip

\bib{La}Laudal, O.~A., ``Formal Moduli of Algebraic Structures,'' 
{\it Lect. Notes 754}, Springer (1979). 

\medskip

\bib{NR}Nijenhuis, A., Richardson, R.W., ``Cohomology and 
deformations in graded Lie algebras,'' {\it Bull. Amer. Math.
Soc.}, {\bf 72} (1966), 1--29.

\medskip

\bib{P}Palamodov, V.P., ``Deformations of complex spaces,'' {\it Russian 
Math. Surveys}, {\bf 31} (1976).

\medskip

\bib{Sch}Schlessinger, M., ``Functors of Artin rings,'' {\it Trans. Amer.
Math. Soc.} {\bf 130} (1968), 208--222.

\bigskip
\parindent=0pt
\obeylines {\sc E-mail}
\quad Alice Fialowski: {\tt fialowsk@cs.elte.hu}
\quad Dmitry Fuchs: {\tt fuchs@math.ucdavis.edu}
\bye